\documentclass[12pt]{article}
\usepackage{latexsym,amsmath,amsfonts,amssymb,amsthm,mathrsfs}
\usepackage[usenames,dvipsnames]{color}
\usepackage{verbatim}

\scrollmode

\newtheorem{thm}{Theorem}[section]
\newtheorem{lem}[thm]{Lemma}
\newtheorem{cor}[thm]{Corollary}
\newtheorem{rem}[thm]{Remark}
\newtheorem{rems}[thm]{Remarks}
\newtheorem{prop}[thm]{Proposition}
\newtheorem{df}[thm]{Definition}

\newtheorem{ex}[thm]{Example}

\newtheorem{ntr}[thm]{Notations \& Remarks}

\newcommand{\wh}{\widehat}
\newcommand{\wt}{\widetilde}
\newcommand{\vS}{\varSigma}
\newcommand{\vO}{\varOmega}
\newcommand{\vT}{\varTheta}
\newcommand{\sq}{\subseteq}

\newcommand{\mf}{\mathfrak}
\newcommand{\vY}{\varUpsilon}
\newcommand{\vX}{\varXi}
\newcommand{\uph}{\upharpoonright}

\def\N{{\mathbb N}}
\def\R{{\mathbb R}}

\def\Q{{\mathbb Q}}

\newcommand{\cnpw}{counting process }

\newcommand{\agpw}{aggregate process }

\newcommand{\cspw}{size process }

\newcommand{\rp}{(\cnp,\csp)}
\newcommand{\rpw}{risk process }
\newcommand{\rpB}{risk process}

\newcommand{\E}{\mathbb{E}}
\newcommand{\M}{\mathbb{M}}
\newcommand{\bdot}{\bullet}
\newcommand{\leb}{\lambda}

\newcommand{\haf}{{\mathcal{F}}}
\newcommand{\xif}{{\xi}}                                                      

\newcommand{\sagp}{S}
\newcommand{\cnp}{N}
\newcommand{\clap}{T}
\newcommand{\agp}{S}
\newcommand{\csp}{X}
\newcommand{\cip}{W}

\newcommand{\Zp}{Z}
\newcommand{\Vp}{V}

\newcommand{\margn}[1]{\marginpar{\texttt{#1}}}

\newcommand{\adc}[1]{\textcolor{blue}{#1}}
\newcommand{\arch}[1]{\textcolor[rgb]{0.5,0,0.5}{#1}}

\newcommand{\del}[1]{\textcolor{red}{#1}}

\newcommand{\mcol}[1]{\textcolor[rgb]{0,0,0}{#1}}

\newcommand{\ins}{\mcol{O}} 
\newcommand{\poins}{\mcol{L_{*}}} 

\newcommand{\lpns}{\mcol{G^{\prime}}} 
\newcommand{\ldns}{\mcol{G^{\prime\prime}}} 
\newcommand{\lrns}{\mcol{G_{*}}} 
\newcommand{\lcns}{\mcol{C_{*}}} 
\newcommand{\pcns}{\mcol{V}} 
\newcommand{\qans}{\mcol{V_{*}}} 
\newcommand{\lcnsq}{\mcol{{\wt{C}}_{*}}} 
\newcommand{\ccns}{\mcol{G_{\mathrm{a}}}} 
\newcommand{\lpnsq}{\mcol{G_{\mathrm{b}}}} 
\newcommand{\ldnsq}{\mcol{\wt{G}^{\prime\prime}}} 
\newcommand{\newd}{{G_{\mathrm{d}}}}
\newcommand{\gplians}{\mcol{{{O}}^{(*)}}} 
\newcommand{\pliani}{\mcol{G_{\mathrm{d}}}} 
\newcommand{\qliani}{\mcol{\wt{O}^{(*)}}} 
\newcommand{\caggns}{\mcol{G_{\mathrm{c}}}}
\newcommand{\machns}{\mcol{\wt{O}_{n,B_n}^{(*)}}} 
\newcommand{\umachns}{\mcol{{O}_{n,B_n}^{(*)}}} 
\newcommand{\umachnsm}{\mcol{{O}_{n,B_n,\gamma,X_1,\ldots,X_n}^{(*)}}}
\newcommand{\qcns}{\mcol{\wt{V}}} 
 
\newcommand{\moic}{\mcol{{\wt{C}}^{(*)}}}
\newcommand{\idnso}{\mcol{L}} 

\newcommand{\T}{{\mathbb{T}}}

\newcommand{\auth}[1]{\textcolor{black}{\texttt{#1}}}
\newcommand{\artp}[1]{\textrm{#1}}
\newcommand{\artn}[1]{\textbf{#1}}
\newcommand{\artj}[1]{\textsl{#1}}
\newcommand{\artb}[1]{\emph{#1}}

\title{{\LARGE\bf A characterization of martingale-equivalent compound mixed Poisson process}}

\author{
{\auth{D.P. Lyberopoulos}\thanks{The author is indebted to the Public Benefit Foundation \textsc{Alexander S. Onassis}, 
which supported this research, under the Programme of Scholarships for Hellenes.}\;\;{\normalsize and}\;\;\auth{N.D. Macheras}}
}

\date{{\footnotesize\textsf{\today}}}
\begin{document}
\maketitle


\begin{abstract}
If a given aggregate process $\agp$ is a compound mixed Poisson process under a probability measure $P$, 
a characterization of all probability measures $Q$ on the domain of $P$, such that $P$ and $Q$ are progressively 
equivalent and $\agp$ remains a compound mixed Poisson process with improved properties, 
is provided. This result generalizes earlier work of  Delbaen \& Haezendonck (1989). 
Implications related to the computation of premium calculation principles in an insurance 
market possessing the property of no free lunch with vanishing risk are also discussed.
\smallskip

\par\noindent{\bf MSC 2010:} Primary 60G55, 91B30 ; secondary 28A50, 60G57, 60A10, 60G44.


\par\noindent{\bf{Key Words}:} {\rm regular conditional probability, compound mixed Poisson process, martingale, martingale-equivalent measures, premium calculation principle}.
\end{abstract}

\section*{Introduction}\label{lm30}

In \cite{dh}, Proposition 2.2, Delbaen \& Haezendonck provided a positive answer to the problem of 
characterizing all martingale-equivalent probability measures $Q$ such that a compound Poisson process under 
an original probability measure $P$ ($P$-CPP for short) will remain a CPP under $Q$. 
This problem raised naturally as the above authors tried ``to create a mathematical framework to deal with finance related to risk processes'' in the 
context of classical Risk Theory. For details see \cite{dh}, Section 1.

In fact, \cite{dh} played a fundamental role in comprehending the interplay between financial and actuarial pricing of insurance
(see e.g. \cite{em00}) and influenced the work of many researchers; we refer for instance to Meister \cite{mei}, Embrechts \& Meister \cite{em97} 
and Holtan \cite{ho}, where the pricing of catastrophe insurance futures and non-life insurance contracts, respectively, was investigated. 
In M{\o}ller \cite{mol}, the ordering of various martingale measures, 
including those suggested in \cite{dh},
was studied in a dynamic reinsurance markets setting. The results of Delbaen \& Haezendonck \cite{dh} were exploited by
Boogaert \& de Waegenaere \cite{bdw} for simulating ruin probabilities, 
and more recently by Yu et al. \cite{yual} to address a problem of capital risk allocation for a redevelopment project.

The purpose of this paper is to provide under weaker 
assumptions than those made in \cite{mei} (see Remark \ref{qr}, (b)) 
an extention of the characterization of Delbaen \& Haezendonck \cite{dh} for the case of compound mixed Poisson processes (CMPPs for short), that is, in the frame of non classical Risk Theory. 
Since conditioning is involved in the definition of CMPPs, it is natural to expect that 
{\rm regular conditional probabilities} will play a key role towards this direction. 
For this reason, we first give a characterization of CMPPs in terms of regular conditional probabilities, see Proposition \ref{12}, via which CMPPs are reduced to ordinary CPPs for the disintegrating measures. 

In order to investigate the existence of progressively equivalent martingale measures $Q$ (see Section \ref{cmppm} for the definition) one has to be able to characterize the Radon-Nikod\'{y}m derivatives of $Q$ with respect to the original measure $P$. This is done in Section \ref{lm34}, Proposition \ref{tsp}.
The latter result together with Proposition \ref{conb}, an existence result for CMPPs, is required in order to  
characterize, in terms of regular conditional probabilities, 
all measures $Q$ such that they are progressively equivalent to an original 
measure $P$, and such that a $P$-CMPP will remain a CMPP under $Q$, see Theorem \ref{qQ}. 
In this way, the main result, Proposition 2.2, of \cite{dh}, becomes a special case of Theorem \ref{qQ}, see Remarks \ref{qr}.

In Section \ref{cmppm}, applying Theorem \ref{qQ}, 
we find out a wide class {\em of canonical} stochastic processes inducing a corresponding one of {\em progressively equivalent martingale measures} (PEMMs for short), see Proposition \ref{cqr}, and possessing the property of {\em no free lunch with vanishing risk} ((NFLVR) for short), see Theorem \ref{qm00}. In particular, in Proposition \ref{cqr} a characterization of all PEMMs $Q$ such that a CMPP under $P$ remains a CMPP under $Q$ is given.

Finally, implications of Theorems \ref{qQ} and \ref{qm00} related to the computation of of premium calculation principles in an insurance 
market possessing the property of (NFLVR) are discussed 
in Section 6.

\section{Preliminaries}\label{lm31}
$\N$ and $\R$ stand for the natural and the real numbers, respectively, while $\R_+:=\{x\in\R:x\geq0\}$. If $d\in\N$, then $\R^d$ denotes the Euclidean space of dimension $d$.

Given a probability space $(\vO,\vS,P)$, a set $N\in\vS$ with $P(N)=0$ is called a $P$-{\bf null set} (or a null set for simplicity). 
For random variables $X,Y:\vO\longrightarrow\R$ we write $X=Y\;$ $P$-almost surely ($P$-a.s. for short), if 
$P(X\neq Y)=0$.

If $A\subseteq\vO$, then $A^c:=\vO\setminus A$, while $\chi_A$ denotes the indicator (or characteristic) function of the set $A$. 
For a map $f:D\longrightarrow{\R}$ and for a non-empty set $A\subseteq{D}$ we denote by $f\upharpoonright{A}$ the restriction of $f$ to $A$. 
The identity map from $\vO$ onto itself is denoted by $id_{\vO}$. 
The $\sigma$-algebra generated by a family $\mathcal{G}$ of subsets of $\vO$ 
is denoted by $\sigma(\mathcal{G})$.

For any Hausdorff topology $\mf{T}$ on $\vO$, by ${\mf B}(\vO)$ is denoted the {\bf Borel $\sigma$-algebra} on $\vO$, i.e. the $\sigma$-algebra generated by $\mf{T}$, while by $\mf{B}:=\mf{B}(\R)$ is denoted the Borel $\sigma$-algebra of subsets of $\R$. 
By $\mathcal{L}^{\ell}(P)$ will be denoted the space of all $\vS$-measurable real-valued functions $f$ on $\vO$ such that $\int|f|^{\ell}dP<\infty$ 
for $\ell\in\{1,2\}$.

Functions that are $P$-$\mbox{a.s.}$ equal are not identified. 
We write $\E_P[X\mid{\mathcal{G}}]$ for a version of a conditional expectation (under $P$) 
of $X\in\mathcal{L}^{1}(P)$ given a $\sigma$-subalgebra ${\mathcal{G}}$ of $\vS$.

Given two probability spaces $(\vO,\vS,P)$ and $(\vY,T,Q)$ as well as a $\vS$-$T$-measurable map $X:\vO\longrightarrow\vY$ we denote by $\sigma(X):=\{X^{-1}(B): B\in T\}$ the $\sigma$-algebra generated by $X$, while 
$\sigma(\{X_i\}_{i\in I}):=\sigma\bigl(\bigcup_{i\in{I}}\sigma(X_i)\bigr)$ stands for the $\sigma$-algebra generated by a family $\{X_i\}_{i\in I}$ of $\vS$-$T$-measurable maps from $\vO$ into $\vY$.
For any given $\vS$-$T$-measurable map $X$ from $\vO$ into $\vY$ denote by $P_X: T\longrightarrow\R$ the image measure of $P$ under $X$. 
By $\mathbf{K}(\theta)$ is denoted an arbitrary probability distribution on $\mf{B}$ with parameter $\theta\in\vX$. In particular, 
$\mathbf{P}(\theta)$ and $\mathbf{Exp}(\theta)$, where $\theta$ is a positive parameter, stand for the law of Poisson and exponential 
distribution, respectively (cf. e.g. \cite{sch}).

Given two real-valued random variable $X,\vT$ on $\vO$,  
a {\bf conditional distribution of $X$ over $\vT$} is a $\sigma(\vT)$-$\mathfrak{B}$-Markov kernel (see \cite{ba}, Definition 36.1 for the definition)
denoted by $P_{X\mid\vT}:=P_{X\mid\sigma(\vT)}$ and satisfying for each $B\in\mf{B}$ the equality
$P_{X\mid\vT}(\bullet,B)=P(X^{-1}(B)\mid\sigma(\vT))(\bullet)$ ${P}\uph\sigma(\vT)$-a.s.. Clearly, for every $\mathfrak{B}_d$-$\mathfrak{B}$-Markov kernel $k$, the map $K(\vT)$ from $\vO\times\mathfrak{B}$ into $[0,1]$ defined by means of
$$
K(\vT)(\omega,B):=(k(\bullet,B)\circ\vT)(\omega)
\quad\mbox{for any}\;\;(\omega,B)\in \vO\times\mathfrak{B}
$$
is a $\sigma(\vT)$-$\mathfrak{B}$-Markov kernel. Then for $\theta=\vT(\omega)$ with $\omega\in\vO$ the probability measures $k(\theta,\bullet)$ are distributions on $\mathfrak{B}$ and so we may write $\mathbf{K}(\theta)(\bullet)$ instead of $k(\theta,\bullet)$. Consequently, in this case $K(\vT)$ will be denoted by 
$\mathbf{K}(\vT)$.

For any real-valued random variables $X$, $Y$ on $\vO$ we say that $P_{X\mid\vT}$ and $P_{Y\mid\vT}$ are $P\uph\sigma(\vT)$-equivalent and we write 
$P_{X\mid\vT}=P_{Y\mid\vT}$ $P\uph\sigma(\vT)$-a.s., if there exists a $P$-null set $N\in\sigma(\vT)$ such that for any $\omega\notin N$ and 
$B\in\mf{B}$ the equality $P_{X\mid\vT}(B,\omega)=P_{Y\mid\vT}(B,\omega)$ holds true.

A family $\{{X}_i\}_{i\in I}$ of random variables 
is {\bf $P$-conditionally identically distributed} over a 
random variable $\vT$, if $P(F\cap{X}_i^{-1}(B))=P(F\cap{X}_j^{-1}(B))$ whenever $i,j\in I$, $F\in\sigma(\vT)$
and $B\in\mf{B}$. 
Furthermore, we say that 
$\{{X}_i\}_{i\in I}$ is {\bf $P$-conditionally (stochastically) independent  given}  
$\vT$, if it is conditionally independent given the $\sigma$-algebra $\sigma(\vT)$; 
for the definition of conditional independence see  e.g. \cite{ct}, page 220.

{\em From now on let $(\vO,\vS,P)$ be an arbitrary but fixed probability space. Unless it is stated otherwise, $\vT$ is a random variable on $\vO$ such that $P_{\vT}\bigl((0,\infty)\bigr)=1$, and we simply write 
``conditionally'' in the place of ``conditionally given $\vT$'' whenever conditioning refers to $\vT$}.


\section{A characterization of compound mixed Poisson processes}\label{lm33}

Let $\cnp:=\{N_t\}_{t\in\R_+}$ be a {\bf ($P$-) \cnpw} with 
exceptional ($P$-) null set $\vO_N$ (cf. e.g. \cite{sch}, page 17 for the definition). Without loss of generality we may and do assume that $\vO_N=\emptyset$. Denote by $\clap:=\{T_n\}_{n\in\N_0}$ and $\cip:=\{W_n\}_{n\in\N}$ the {\bf arrival process} and 
{\bf interarrival process}, respectively, associated with $\cnp$ (cf. e.g. \cite{sch}, page 6 for the definitions).
Also let $\csp:=\{X_n\}_{n\in\N}$ be the {\bf (size process} with all $X_n$ positive, and $\agp:=\{S_t\}_{t\in\R_+}$ the {\bf aggregate  process} induced by the counting process $\cnp$ and the size process $\csp$ (cf. e.g. \cite{sch}, page 103 for the definitions). 
For the definition of a risk process $\rp$ on $(\vO,\vS,P)$ we refer to \cite{sch}, page 127.

The \cnpw $\cnp$ is said to be a {\bf mixed Poisson process} on $(\vO,\vS,P)$ with parameter $\vT$ (or a $P$-MPP($\vT$) for short),  if it has conditionally stationary independent increments, such that
$$
P_{N_t\mid\vT}=\mathbf{P}(t\vT)\quad{P\upharpoonright\sigma(\vT)-\mbox{a.s.}}
$$ 
holds true for each $t\in(0,\infty)$.
\smallskip

In particular, if the distribution of $\vT$ is degenerate at $\theta_0>0$ (i.e. $P_{\vT}(\{\theta_0\})=1$), then $\cnp$ is a {\bf $P$-Poisson process} with parameter $\theta_0$ (or a $P$-PP($\theta_0$) for short).
\smallskip

An \agpw $\agp$ is said to be a {\bf compound mixed Poisson process} on $(\vO,\vS,P)$ with parameters $\vT$ and $P_{X_1}$ (or else a $P$-CMPP$(\vT,P_{X_1})$ for short), if it is induced by a $P$-risk process $\rp$ such that $\cnp$ is a $P$-MPP($\vT$) 
(cf. e.g. \cite{sch}, Section 5.1). 

In particular, if the distribution of $\vT$ is degenerate at $\theta_0>0$ then $\agp$ is said to be a {\bf compound Poisson process} on $(\vO,\vS,P)$ with parameters $\theta_0$ and $P_{X_1}$ (or else a $P$-CPP$(\theta_0,P_{X_1})$ for short).\medskip

The following conditions will serve as useful assumptions in the study of CMPPs: 
\begin{description}
\item[(a1)] The processes $\cip$ and $\csp$ are $P$-conditionally mutually independent.
\item[(a2)] The random variable $\vT$ and the sequence $\csp$ are $P$-(unconditionally) independent.
\end{description}
{\em Next, whenever condition 
\textnormal{(a1) and (a2)} 
holds true we shall write that the quadruplet $(P,\cip,\csp,\vT)$ 
or (if no confusion arises) the probability measure $P$ satisfies
\textnormal{(a1)} and \textnormal{(a2)}, respectively}.

Consider now a second arbitrary but fixed probability space $ (\vY,T,Q)$. The following definition is a special instance of that in \cite{fr4}, 452E, proper for our investigation.

\begin{df}\label{dis} 
\normalfont
A \textbf{regular conditional probability} (r.c.p. for short) \textbf{of $P$ over $Q$} is a family $\{P_{y}\}_{y\in\vY}$ of probability measures $P_y:\vS\longrightarrow\R$ such that 
\begin{description}
\item[(d1)] for each $D\in\vS$ the function $P_{\cdot}(D):\vY\longrightarrow\R$ is $T$-measurable;
\item[(d2)] $\int P_{y}(D)Q(dy)=P(D)$ for each $D\in\vS$.
\end{description}
We could use the name of {\em disintegration} instead, 
but it seems that it is better to reserve that term to the general case when
$P_y$'s may be defined on different domains (see \cite{pa}).

If $f:\vO\longrightarrow\vY$ is an inverse-measure-preserving map (i.e. $P(f^{-1}(B))=Q(B)$ for each $B\in{T}$), a r.c.p. $\{P_{y}\}_{y\in\vY}$ of $P$ over $Q$ is called \textbf{consistent} with $f$ if, for each $B\in{T}$, the equality $P_{y}(f^{-1}(B))=1$ 
holds for $Q$-almost all ($Q$-a.a. for short) $y\in B$.

We say that a r.c.p. $\{P_{y}\}_{y\in\vY}$ of $P$ over $Q$ consistent with $f$ 
is \textbf{essential unique}, if for any other r.c.p. $\{P_{y}^{\prime}\}_{y\in\vY}$ of $P$ over $Q$ consistent with $f$ there exists a 
\mcol{$P_{\vT}$-null set $N\in{T}$ such that for any $\theta\notin N$ the equality $P_y=P_y^{\prime}$ holds true.} 
\end{df}

\begin{rem}\label{magd} 
\normalfont
If $\vS$ is countably generated and 
$(\vO,\vS,P)$ or $P$ is perfect (see \cite{fa}, page 291 for the definition),
then there always exists a r.c.p. $\{P_{y}\}_{y\in\vY}$ of $P$ over $Q$ consistent with any inverse-measure-preserving map $f$ from $\vO$ into $\vY$ providing that $T$ is countably generated (see \cite{fa}, Theorems 6 and 3). Note that the most important applications in Probability Theory are still rooted in the case of standard Borel spaces $(\vO,\vS)$, that is, of spaces being isomorphic to $(Z,\mf{B}(Z))$, where $Z$ is some Polish space; hence of spaces satisfying always the above mentioned assumptions concerning $P$, $\vS$ and $T$. It is well-known that any Polish space is standard Borel; in particular, $\R^{d}$ and $\R^{\N}$ are such spaces. If $(\vO,\vS)$ and $(\vY,T)$ are non-empty standard Borel spaces, then there always exists an essentially unique r.c.p. $\{P_{y}\}_{y\in\vY}$ of $P$ over $Q$ consistent with any inverse-measure-preserving map $f$ from $\vO$ into $\vY$ (cf. e.g. \cite{fr4}, 452X(m)).
\end{rem}

If $\cnp$ is a $P$-MPP$(\vT)$, we then get by e.g. \cite{sch}, Lemma 4.2.5, that the {\bf explosion} $E:=\{\sup_{n\in\N}T_n<\infty\}$ is a $P$-null set. 
{\em Without loss of generality we may and do consider explosion equal to the empty set.} 

{\em Throughout what follows we put $\vY:=(0,\infty)$ and assume that there exists a r.c.p. $\{P_{\theta}\}_{\theta\in\vY}$ of $P$ over $P_{\vT}$ consistent with $\vT$}.

\begin{lem}\label{11}
\begin{enumerate}
\item
The following conditions are equivalent:
\begin{itemize}
\item[(a)] Condition \textnormal{(a1)};
\item[(b)] $\cnp$ and $\csp$ are $P$-conditionally mutually independent;
\item[(c)] 
there exists a $P_{\vT}$-null set ${\lpns}\in\mf{B}(\vY)$ such that for any $\theta\notin{\lpns}$ the processes $\cnp$ and $\csp$ are 
$P_{\theta}$-mutually independent.
\end{itemize}
\item
Condition \textnormal{(a2)} implies that the process $\csp$ is $P$-i.i.d. if and only 
if it is $P$-conditionally i.i.d. if and only if there exists a $P_{\vT}$-null set ${\ldns}\in\mf{B}(\vY)$ such that 
for any $\theta\notin{\ldns}$ the process $\csp$ is $P_{\theta}$-i.i.d.. 
\item
Conditions \textnormal{(a1)} and \textnormal{(a2)} imply that the pair $\rp$ is a $P$-\rpw if and only if 
there exists a $P_{\vT}$-null set ${\lrns}\in\mf{B}(\vY)$ such that for any $\theta\notin{\lrns}$ the pair $\rp$ 
is a $P_{\theta}$-risk process.
\end{enumerate}
\end{lem}

{\bf Proof.} 
Ad $(i)$: First note that $\sigma(\mathcal{F}^T_n)=\sigma(\mathcal{F}^W_n)$ 
for each $n\in\N$ (cf. \cite{sch}, Lemma 1.1.1), which implies that $\sigma(\cip)=\sigma(\clap)=\sigma(\cnp)$, 
where the last equality is an immediate consequence of \cite{sch}, Lemma 2.1.3; hence implication $(b)\Longrightarrow(a)$ follows. The inverse 
implication is immediate since 
\begin{equation}\label{Z0}
N_t
=\sum_{n=0}^{\infty}n\chi_{\{T_n\leq t<T_{n+1}\}}
=\sum_{n=0}^{\infty}n\chi_{\{\sum_{k=1}^nW_k\leq t<\sum_{k=1}^{n+1}W_k\}}
\end{equation}
for each $t\in\R_+$ (cf. e.g. \cite{sch}, Theorem 2.1.1 and Lemma 2.1.2). 
The equivalence $(a)\Longleftrightarrow(c)$ follows by 
\cite{lm1}, Lemma 4.1 together with \cite{lm1err}
\smallskip 

Ad $(ii)$: 
Assume that \textnormal{(a2)} holds true. Then $\csp$ is $P$-independent if and only if it is $P$-conditionally independent. 
But applying a monotone class argument we easily conclude that the latter is equivalent to the fact that
there exists a $P_{\vT}$-null set ${\ins}\mcol{:={\ins}_{\N}}\in\mf{B}(\vY)$ such that for any $\theta\notin{\ins}$ the process $\csp$ is $P_{\theta}$-independent.

By (a2) the fact that all random variables $X_n$ are $P$-identically distributed is equivalent to the fact that they are $P$-conditionally identically distributed, which again by a monotone class argument  
equivalently yields that there exists a $P_{\vT}$-null set $\mcol{{L}:={L}_{\N}}\in\mf{B}(\vY)$ such that for any $\theta\notin{L}$ the process $\csp$ is $P_{\theta}$-identically distributed. So, putting ${\ldns}:={\ins}\cup{L}$ assertion $(ii)$ follows. 
\smallskip

Ad $(iii)$: 
Assume that $P$ satisfies \textnormal{(a1)} and \textnormal{(a2)} as well as that there exists a $P_{\vT}$-null set ${\lrns}\in\mf{B}(\vY)$ such that for any $\theta\notin{\lrns}$ the pair $\rp$ is a $P_{\theta}$-risk process. Then by $(ii)$ we get that 
$\csp$ is a sequence of positive $P_{\theta}$-i.i.d. random variables on $\vO$ 
for any $\theta\notin{\lrns}$ if and only if it does so on $(\vO,\vS,P)$. Also note that according to $(i)$ the 
fact, that the processes $\csp$ and $\cnp$ 
are mutually independent under  $P_{\theta}$ for any $\theta\notin{\lrns}$, 
is equivalent to condition \textnormal{(a1)}, which together with condition \textnormal{(a2)} yields that 
the processes $\csp$ and $\cnp$ are mutually independent under $P$.
Consequently, the pair $\rp$ is a $P$-\rpB. 

The inverse implication follows for any $\theta\notin{\lrns}:={\lpns}\cup{\ldns}$ by assertions $(i)$ and $(ii)$, since 
clearly $\cnp$ is a $P_{\theta}$-counting process for any $\theta\notin{L}$.\hfill$\Box$

\begin{prop}\label{12}
Assume that $P$ satisfies conditions \textnormal{(a1)} and \textnormal{(a2)}.
Then the \agpw $\agp$ is a $P$-CMPP$(\vT,P_{X_1})$ if and only if 
there exists a $P_{\vT}$-null set ${\lcns}\in\mf{B}(\vY)$ such that for any 
$\theta\notin{\lcns}$ the family $\agp$ is a $P_{\theta}$-CPP$(\theta,(P_{\theta})_{X_1})$. 
\end{prop}

{\bf Proof.} Assume that $\agp$ is a $P$-CMPP$(\vT,P_{X_1})$, which is equivalent to the fact that the pair $\rp$ is a $P$-\rpw and $\cnp$ is a $P$-MPP with parameter $\vT$. According to Lemma \ref{11}, $(iii)$ and 
\cite{lm1}, Proposition 4.4, the latter equivalently yields that there exist two $P_{\vT}$-null sets ${\lrns}$ 
and ${\poins}$ in $\mf{B}(\vY)$ such that for any $\theta\notin{\lrns}$ the pair $\rp$ is a $P_{\theta}$-\rpB, and 
for any $\theta\notin{\poins}$ the family $\cnp$ is a $P_{\theta}$-PP($\theta$), respectively. 
So, \mcol{if we let} ${\lcns}={\lrns}\cup{\poins}$ we then equivalently get that $\agp$ is a $P_{\theta}$-CPP$(\theta,(P_{\theta})_{X_1})$ 
for any $\theta\notin{\lcns}$.\hfill$\Box$

\section{Change of measures for compound mixed Poisson processes}\label{lm34}

The main result of this section, Proposition \ref{tsp}, allows us to explicitly calculate Radon-Nikod\'{y}m derivatives for the most important insurance risk processes.

Let $\T\subseteq\R_+$ with $0\in\T$ and let $\ell\in\{1,2\}$. 
For a process $Z_{\T}:=\{Z_t\}_{t\in\T}$ denote by $\mathcal{F}^Z_{\T}:=\{\mathcal{F}^Z_t\}_{t\in\T}$
the canonical filtration of $Z_{\T}$. For $\T=\R_+$ 
write $Z$ and $\mathcal{F}^Z$ in the place of $Z_{\R_+}$ and $\mathcal{F}^Y_{\R_+}$, respectively.
Write also $\mathcal{F}:=\{\mathcal{F}_t\}_{t\in\R_+}$, where $\mathcal{F}_t:=\sigma\bigl(\mathcal{F}^{\sagp}_t\cup\sigma(\vT)\bigr)$ for the canonical filtration of $\sagp$ and $\vT$, $\mathcal{F}_{\infty}^{\sagp}:=\sigma(\mathcal{F}^{\sagp})$ and 
$\mathcal{F}_{\infty}:=\sigma\bigl(\mathcal{F}^{\sagp}_{\infty}\cup\sigma(\vT)\bigr)$ for simplicity.
Recall that a {\bf martingale in} $\mathcal{L}^{\ell}(P)$ {\bf adapted to the filtration} 
$\mathcal{Z}_{\T}$, or else a $\mathcal{Z}_{\T}$-{\bf martingale} in $\mathcal{L}^{\ell}(P)$, is a process $Z_{\T}:=\{Z_t\}_{t\in\T}$ of  
real-valued random variables in $\mathcal{L}^{\ell}(P)$ such that $Z_t$ is $\mathcal{Z}_t$-measurable for each $t\in\T$ and whenever $s\leq{t}$ in $\T$ and $E\in\mathcal{Z}_s$ then $\int_EZ_sdP=\int_EZ_tdP$.
The latter condition is called the {\bf martingale property} (cf. e.g. \cite{sch}, page 25). 
For $\mathcal{Z}=\mathcal{F}$ we simply say that $Z$ is a
martingale in $\mathcal{L}^{\ell}(P)$.

\begin{rem}\label{qS1} 
\normalfont
\noindent
For any $n\in\N$ the random variable $X_n$ is 
$\haf_{T_n}^{\sagp}$-measurable, where 
$$
\haf_{T_n}^{\sagp}:=\{A\in\vS: A\cap\{T_n\leq{t}\}\in\haf_t^{\sagp}\;\;\mbox{for every}\;\; t\in\R_+\},
$$
and for any $t\in\R_+$ the random variable $X_{N_t}$ is 
$\haf_t^{\sagp}$-measurable.

In fact, it follows by \cite{sch}, Lemma 2.1.2 
that all random variables $T_n$ are 
$\haf^{\sagp}$-stopping times.
Furthermore, $\agp$ is right-continuous, since $\cnp$ is so. The latter together with the fact that $T_{n-1}<T_n$ for any $n\in\N$ yields that the random variables $S_{T_n}$ and $S_{T_{n-1}}$ are 
$\haf_{T_n}^{\sagp}$-measurable 
for each $n\in\N$ (cf. e.g. \cite{ks}, Chapter 1, Propositions 2.18, 1.13 and Lemma 2.15). 
Thus, taking into account that 
$X_n=S_{T_n}-S_{T_{n-1}}$ since $N_{T_n}=n$ for each $n\in\N$, we deduce that $X_n$ is 
$\haf_{T_n}^{\sagp}$-measurable for any $n\in\N$. 

But for all $n\in\N_0$ and $t\in\R_+$ we have 
$\{N_t=n\}=\{T_n\leq t<T_{n+1}\}\in\haf_t^{\sagp}$ 
(see \cite{sch}, Lemma 2.1.2 for the equality), 
implying that $X_{N_t}^{-1}(B)\cap\{N_t=n\}\in\haf_t^{\sagp}$ 
for each $B\in\mf{B}(\vY)$ (see \cite{ks}, Chapter 1, Lemma 2.15). 
Consequently, the 
$\haf_t^{\sagp}$-measurability of each random variable $X_{N_t}$ follows.
\end{rem}
 
\begin{df}\label{emm}
\normalfont
Let 
$P$, $Q$ be two probability measures on $\vS$ and 
$\{Y_t\}_{t\in{\T}}$ a process on $(\vO,\vS)$. Then
$P$ and $Q$ are said to be {\bf progressively equivalent} 
if $P$ and $Q$ are equivalent (in the sense of absolute continuity) 
on each $\mathcal{Z}_t$, in symbols $Q\stackrel{pr}{\sim}P$. If $P$ and $Q$ are equivalent on $\vS$ we write $P\sim{Q}$.
\end{df}

\begin{ntr}\label{md}
\normalfont
\textbf{(a)}
The class of all $\mf{B}(\vY\times\vY)$-measurable real-valued 
functions $\beta$ on $\vY\times\vY$ such that $\beta(x,\theta):=\alpha(\theta)+\gamma(x)$ for each $x,\theta\in\vY$, where $\alpha$, $\gamma$ are 
$\mf{B}(\vY)$-measurable functions from $\vY$ into $\R$ with 
$\E_P[e^{\gamma(X_1)}]=1$, is denoted by $\mathcal{F}_P:=\mathcal{F}_{P,X_1,\vT}$. 
The class of all real-valued functions $\xif$ on $\vY$ such that $P_{\vT}(\{\xif>0\})=1$   
and $\E_P[\xif(\vT)]=1$ is denoted by $\mathcal{R}_+(\vY):=\mathcal{R}_+(\vY,\mathfrak{B}(\vY), P_{\vT})$.  

\textbf{(b)}
Denote by $\mathfrak{M}_+(\vY):=\mathfrak{M}_+(\vY,\mathfrak{B}(\vY))$ the class of all positive, $\mathfrak{B}(\vY)$-measurable 
functions on $\vY$.  
Then for each $g\in\mf{M}_+(\vY)$ the class of all probability measures $Q$ on $\vS$ satisfying \textnormal{(a1)} and \textnormal{(a2)}, 
such that $Q\stackrel{pr}{\sim}P$ and 
$\agp$ is a $Q$-CMPP$(g(\vT),Q_{X_1})$ is denoted by 
$\mathcal{M}_{\sagp,g}:=\mathcal{M}_{\sagp,g,P,{X_1},\vT}$. 
For $g:=id_{\vY}$ put $\mathcal{M}_{\sagp}:=\mathcal{M}_{\sagp,id_{\vY}}$.

\textbf{(c)}
For any $\theta\in\vY$ denote by $\wt{\mathcal{M}}_{\sagp,g}:=\wt{\mathcal{M}}_{\sagp,g}(\theta)$ 
the class of all probability measures 
$Q_\theta$ on $\vS$ such that $Q_{\theta}\stackrel{pr}{\sim} P_{\theta}$ and $\sagp$ is a 
$Q_{\theta}$-CPP$(g(\theta),(Q_{\theta})_{X_1})$. For any given $Q\in\mathcal{M}_{\sagp,g}$, if $\{Q_{\theta}\}_{\theta\in\vY}$ is a r.c.p. of $Q$ 
over $Q_{\vT}$ consistent with $\vT$, it 
follows by Proposition \ref{12} that there exists a $P_{\vT}$-null set $\wh{L}\in\mf{B}(\vY)$ such that 
$Q_{\theta}\in{\wt{\mathcal{M}}}_{\sagp,g}$ for any $\theta\notin\wh{L}$. 
\end{ntr}

{\em Henceforth, unless stated otherwise, $P\in\mathcal{M}_S$ is the initial probability measure on $\vS$ 
under which the family $S$ is a CMPP$(\vT,P_{X_1})$}.

The next result provides the one direction of the desired characterization. 

\begin{prop}\label{tsp} 
For given $g\in\mf{M}_+(\vY)$ let 
$Q$ be a probability measure on $\vS$ satisfying conditions \textnormal{(a1)}, \textnormal{(a2)} 
and such that the family $S$ is a $Q$-$CMPP(g(\vT), Q_{X_1})$. Suppose that there exists a r.c.p. $\{Q_{\theta}\}_{\theta\in\vY}$ 
of $Q$ over $Q_{\vT}$ consistent with $\vT$. 
Then the following are equivalent:
\begin{enumerate}
\item
$Q\stackrel{pr}\sim P$.
\item
$Q_{X_1}\sim P_{X_1}$ and $Q_{\vT}\sim P_{\vT}$.
\item
There exist an essentially unique $\beta\in\mathcal{F}_P$ such that  
$$
g(\theta)=\theta{e}^{\alpha(\theta)}\quad\mbox{for any}\quad\theta\in\vY,\quad\mbox{and}\quad\gamma=\ln f,\leqno(*)
$$
where $f$ is a \mcol{$P_{X_1}$-a.s. positive} Radon-Nikod\'{y}m derivative of $Q_{X_1}$ with respect to $P_{X_1}$, 
and there exists a $P_{\vT}$-null set ${\qans}\in\mf{B}(\vY)$ such that for any $\theta\notin{\qans}$ 
$$
Q_{\theta}(A)=\int_{A}{\wt{M}}_t^{(\beta)}(\theta)dP_{\theta}\qquad\forall\;0\leq s\leq t\;\;\;\forall\;A\in{\haf}_s,\eqno(M_{\theta})
$$
where 
$\wt{M}_t^{(\beta)}(\theta):=e^{S_t^{(\beta)}(\theta)-t\theta(e^{\alpha(\theta)}-1)}$ with 
$S_t^{(\beta)}(\theta):=\sum_{k=1}^{N_t}\beta(X_k,\theta)$,
and the family 
$\wt{M}^{(\beta)}(\theta):=\{{\wt{M}}_t^{(\beta)}(\theta)\}_{t\in\mathbb{R}_+}$ is a 
martingale
in $\mathcal{L}^1(P_{\theta})$.
\item
There exist an essentially unique 
pair $(\beta,\xif)\in\mathcal{F}_P\times\mathcal{R}_+(\vY)$, where $\xif$ is a Radon-Nikod\'{y}m derivative 
$\xif$ of $Q_{\vT}$ with respect to $P_{\vT}$, such that
$$
g(\vT)=\vT e^{\alpha(\vT)}\quad\mbox{and}\quad\gamma=\ln f,\leqno(**)
$$
where $f$ is a \mcol{$P_{X_1}$-a.s. positive} Radon-Nikod\'{y}m derivative of $Q_{X_1}$ with respect to $P_{X_1}$, 
and such that
$$
Q(A)=\int_{A}{M}_t^{(\beta)}(\vT)dP\qquad\forall\;0\leq s\leq t\;\;\;\forall\;A\in{\haf}_s,
\eqno(M_{\xif})
$$
where 
${{M}}_t^{(\beta)}(\vT):=\xif(\vT){\wt{M}}_t^{(\beta)}(\vT):=\xif(\vT)e^{S_t^{(\beta)}(\vT)-t\vT(e^{\alpha(\vT)}-1)}$, 
and the family $M^{(\beta)}(\vT):=\{M_t^{(\beta)}(\vT)\}_{t\in\mathbb{R}_+}$ is a 
martingale in $\mathcal{L}^1(P)$.
\end{enumerate}
\end{prop}

{\bf Proof.}
Ad $(i)\Longrightarrow (ii)$: 
Since for each $t\in\R_+$ we have $Q\upharpoonright{\haf}_t\sim P\upharpoonright{\haf}_t$ and 
$\haf_t^{\sagp}\subseteq{\haf}_t$, we get $Q\upharpoonright\haf_t^{\sagp}\sim P\upharpoonright\haf_t^{\sagp}$. 

Thus, applying \cite{dh}, Lemma 2.1, we get the first part of assertion $(ii)$. 
The second part of assertion $(ii)$ is an immediate consequence of $Q\stackrel{pr}{\sim}P$ together with 
$\sigma(\vT)\subseteq{\haf}_t$ for each $t\in\R_+$.

Ad $(ii)\Longrightarrow (iii)$: 
\noindent{\bf (a)}
If $\alpha(\theta):=\ln(\frac{g(\theta)}{\theta})$ for each $\theta\in\vY$, then there exists a $P_{\vT}$-null set ${\ccns}\in\mf{B}(\vY)$ such that for any 
$\theta\notin{\ccns}$ and for all $s,t\in\R_+$ with $s\leq{t}$ 
the equality
$$
Q_{\theta}(N_t-N_s=n)=e^{n\alpha(\theta)-(t-s)\theta[e^{\alpha(\theta)}-1]}P_{\theta}(N_t-N_s=n)
$$
holds true for each $n\in\N$.

In fact, since $\sagp$ is both a $P$-CMPP$(\vT,P_{X_1})$ and a $Q$-CMPP$(g(\vT),Q_{X_1})$, it follows by Proposition \ref{12} that there exist a 
$P_{\vT}$- and a $Q_{\vT}$-null set ${\lcns}$ and ${\lcnsq}$ in $\mf{B}(\vY)$ such that for any 
$\theta\notin{\lcns}$ and $\theta\notin{\lcnsq}$ it is a $P_{\theta}$-CPP$(\theta,(P_{\theta})_{X_1})$ and a 
$Q_{\theta}$-CPP$(g(\theta),(Q_{\theta})_{X_1})$, respectively. 
Consequently, putting ${\ccns}:={\lcns}\cup{\lcnsq}$  step (a) follows. 

\noindent{\bf{(b)}}
There exist a $P_{X_1}$-a.s. positive function $f$ such that $Q_{X_1}(\mcol{D})=\int_{\mcol{D}}fdP_{X_1}$ for each $\mcol{D}\in\mf{B}(\vY)$, and 
a $P_{\vT}$-null set $G_b\in\mf{B}(\vY)$ such that for any $\theta\notin{G_b}$, for every $n\in\N$ and 
$B_n\in{\mathcal{F}}_n^{X}$ the equalities
\begin{equation}\label{antid}
Q(B_n)=Q_{\theta}(B_n)
=\E_{P_{\theta}}\bigl[\chi_{B_n}e^{\sum_{j=1}^{n}\gamma(X_j)}\bigr]
=\E_{P}\bigl[\chi_{B_n}e^{\sum_{j=1}^{n}\gamma(X_j)}\bigr],
\end{equation}
hold true with $\gamma:=\ln f$.

In fact, first note that assumption $Q_{X_1}{\sim}P_{X_1}$ implies by the Radon-Nikod\'{y}m Theorem the existence of a $P_{X_1}$-a.s. positive function $f$, which is a Radon-Nikod\'{y}m derivative $f$ of $Q_{X_1}$ with respect to $P_{X_1}$.

We may and do assume without loss of generality that the $P_{X_1}$-null set $\{y\in\vY:f(y)=0\}\in\mf{B}(\vY)$ is the empty set. 

Also note that $\csp$ is $Q$-i.i.d. by assumption. The latter together with Lemma \ref{11}, $(ii)$ implies that there exists a $Q_{\vT}$-null set ${\ldnsq}\in\mf{B}(\vY)$ such that for any $\theta\notin{\ldnsq}$ the sequence $\csp$ is $Q_{\theta}$-i.i.d.; hence the second equality follows in the same way as in \cite{dh}, (2.21).

In order to show the first and last equality of condition (\ref{antid}), 
fix now on $n\in\N$ and $B_n\in{\mathcal{F}_n^{X}}$. 
Since $Q$ satisfies (a2), applying 
\cite{lm1}, Lemma 3.5, we get for \mcol{any} ${D}\in\mf{B}(\vY)$ that 
$$
\int_{{D}}Q_{\theta}(B_n)Q_{\vT}(d\theta)=\int_{\vT^{-1}({D})}Q(B_n\mid\vT)dQ
\stackrel{\mathrm{(a2)}}{=}
\int_{\vT^{-1}({D})}Q(B_n)dQ=\int_{{D}}Q(B_n)Q_{\vT}(d\theta);
$$
hence there exists a $Q_{\vT}$-null set ${\machns}\in\mf{B}(\vY)$ such that for any $\theta\notin{\machns}$ the first equality of (\ref{antid}) holds 
true. But since each $\mathcal{F}_n^{X}$ is countably generated, it can be easily proven by a monotone class argument that there exists a $Q_{\vT}$-null set ${\qliani}\in\mf{B}(\vY)$ such that for any $\theta\notin{\qliani}$ the first equality of (\ref{antid}) holds true.

Since $P$ satisfies (a2), applying \mcol{again \cite{lm1}, Lemma 3.5}, we get for \mcol{any} ${D}\in\mf{B}(\vY)$ that
\[
\int_{\vT^{-1}({D})}\E_{P}\bigl[\chi_{B_n}e^{\sum_{j=1}^{n}\gamma(X_j)}\bigr]dP
=\int_{\vT^{-1}({D})}\E_{P_{\bdot}}
\bigl[\chi_{B_n}e^{\sum_{j=1}^{n}\gamma(X_j)}\bigr]\circ\vT{dP}
\]
or equivalently
$$
\int_{{D}}\E_P\bigl[\chi_{B_n}e^{\sum_{j=1}^{n}\gamma(X_j)}\bigr]{P_{\vT}}(d\theta)
=\int_{{D}}
\E_{P_{\theta}}\bigl[\chi_{B_n}e^{\sum_{j=1}^{n}\gamma(X_j)}\bigr]{P_{\vT}(d\theta)};
$$
hence there exists a $P_{\vT}$-null set ${\umachns}:={\umachnsm}\in\mf{B}(\vY)$ such that for any $\theta\notin{\umachns}$ the last equality of (\ref{antid}) holds true. Again by a monotone class argument it follows that there exists a $P_{\vT}$-null set ${\gplians}\in\mf{B}(\vY)$ such that for any $\theta\notin{\gplians}$ the last equality of (\ref{antid}) holds true. 
So, putting $G_b:={\ldnsq}\cup{\qliani}\cup{\gplians}$, step (b) follows.

\noindent{\bf{(c)}}
Define the set-function $\beta:\vY\times\vY\longrightarrow\mathbb{R}$ by means of $\beta(x,\theta):=\alpha(\theta)+\gamma(x)$ for each $(x,\theta)\in\vY\times\vY$. 
Then $\beta\in{\mathcal{F}_P}$, since by step (b) we have 
$\E_{P_{\theta}}[e^{\gamma(X_1)}]=\E_P[e^{\gamma(X_1)}]=1$ for any $\theta\notin{G_b}$.

{\bf (d)} 
Consider the family of sets 
$$
\mathcal{K}_s^{\sagp}:=\Bigl\{\bigcap_{i=1}^{m} A_i:A_i\in\mathcal{G}_s^{\sagp},\;m\in\N\Bigr\},
$$
where $\mathcal{G}_s^{\sagp}:=\bigcup_{u\leq s}\sigma(S_u)$.

Then there exists
a $P_{\vT}$-null set ${\newd}\in\mf{B}(\vY)$ such that for any 
$\theta\notin{\newd}$ and for any 
$A\in\mathcal{K}_s^{\sagp}$ 
condition $(M_{\theta})$ holds true.

In fact, since $Q$ satisfies \textnormal{(a1)} we get by Lemma \ref{11}, $(i)$ that there exists a $Q_{\vT}$-null set $\wt{G}_d\in\mf{B}(\vY)$ such that for 
any $\theta\notin\wt{G}_d$ the processes $\cnp$ and $\csp$ are $Q_{\theta}$-independent. Let 
$A\in\mathcal{K}_s^{\sagp}$. 

It then follows that there exists a number $m\in\N$ such that $A=\bigcap_{i=1}^mA_i$. 
By induction, it is sufficient to give the proof only for $m=2$.

Then there exist $u_i\in[0,s]$ such that $A_i\in\sigma(S_{u_i})$ 
for $i\in\{1,2\}$ and 
$A=A_1\cap A_2=S_{u_1}^{-1}(B_1)\cap S_{u_2}^{-1}(B_2)$ for some $B_1,B_2\in\mf{B}(\vY)$. 
But since we have assumed that explosion is equal to the empty set, it follows that
$$
A_i=\bigcup_{{m}_i\in\N_0}(\{N_{u_i}={m}_i\}\cap C_{{m}_i}),
$$  
where $C_{{m}_i}:=\{\sum_{j=1}^{{m}_i} X_j\in B_i\}\in{\mathcal{F}_{{m}_i}^{X}}$ for $i\in\{1,2\}$. 

Next without loss of generality we may and do assume that $u_1\leq u_2$. We then get
\begin{eqnarray}
\lefteqn{Q_{\theta}(A\cap\{N_t=n\})}\nonumber\\
&=&
Q_{\theta}\Bigl(\bigcup_{{m}_1=0}^{n}\bigcup_{{m}_2=0}^{n}(C_{{m}_1}\cap\{N_{u_1}={m}_1\}\cap C_{{m}_2}\cap\{N_{u_2}={m}_2\}\cap\{N_t=n\})\Bigr)\nonumber\\
&=&\sum_{{m}_1=0}^{n}\sum_{{m}_2=0}^{n}Q_{\theta}(C_{{m}_1}\cap C_{{m}_2})Q_{\theta}(\{N_{u_1}={m}_1,N_{u_2}={m}_2,N_t=n\}).
\label{at1}
\end{eqnarray}

But by virtue of (a) it follows that for any $\theta\notin{\ccns}$ the \agpw 
$\sagp$ is both a $P_{\theta}$-CPP$(\theta,(P_{\theta})_{X_1})$ and a 
$Q_{\theta}$-CPP$(g(\theta),(Q_{\theta})_{X_1})$, which yields 
\begin{eqnarray}
\varpi_{\theta}&:=&Q_{\theta}(\{N_{u_1}={m}_1,N_{u_2}={m}_2,N_t=n\})\nonumber\\
&=&e^{{m}_1\alpha(\theta)-u_1\theta[e^{\alpha(\theta)}-1]}P_{\theta}(\{N_{u_1}=m_1\})\cdot
e^{({m}_2-{m}_1)\alpha(\theta)-(u_2-u_1)\theta[e^{\alpha(\theta)}-1]}\nonumber\\
& & \cdot P_{\theta}(\{N_{u_2}-N_{u_1}=m_2-m_1\})\cdot e^{({m}_2-{m}_1)\alpha(\theta)-(t-u_2)\theta[e^{\alpha(\theta)}-1]}\nonumber\\
& & \cdot P_{\theta}(\{N_{u_2}-N_{u_1}=m_2-m_1\})\nonumber\\
&=&\E_{P_{\theta}}[\chi_{\{N_{u_1}={m}_1,N_{u_2}={m}_2,N_t=n\}}e^{n\alpha(\theta)-t\theta[e^{\alpha(\theta)}-1]}].\label{at2}
\end{eqnarray}

Furthermore, step (b) yields for any $\theta\notin{G_b}$ that
\begin{equation}
Q_{\theta}(C_{m_1}\cap{C}_{m_2})=\E_{P_{\theta}}[\chi_{C_{m_1}\cap{C}_{m_2}}e^{\sum_{j=1}^{n}\gamma(X_j)}].\label{at3}
\end{equation}

Then equality (\ref{at1}) together with equalities (\ref{at2}),(\ref{at3}) and step (c) yields for any 
$\theta\notin{\newd}:={\ccns}\cup{\lpnsq}\cup\wt{G}_d$ that
\begin{eqnarray*}
Q_{\theta}(A\cap\{N_t=n\})
&=&\E_{P_{\theta}}\bigl[\chi_{\bigcup_{{m}_1=0}^{n}\bigcup_{{m}_2=0}^{n}C_{{m}_1}\cap C_{{m}_2}\cap\{N_{u_1}={m}_1\}\cap\{N_{u_2}={m}_2\}\cap\{N_t=n\}}\\
&\;&\cdot{e}^{\sum_{j=1}^{n}[\alpha(\theta)+\gamma(X_j)]-t\theta[e^{\alpha(\theta)}-1]}\bigr]\\
&=&\E_{P_{\theta}}\bigl[\chi_{A\cap\{N_t=n\}}e^{S_t^{(\beta)}(\theta)-t\theta\E_{P_{\theta}}[e^{\beta^{\theta}(X_1)}-1]}\bigr]
=\E_{P_{\theta}}[\chi_{A\cap\{N_t=n\}}{\wt{M}}_t^{(\beta)}(\theta)],
\end{eqnarray*}
which implies that
$$
Q_{\theta}(A)=\sum_{n\in\N_0} Q_{\theta}(A\cap\{N_t=n\})=\int_A{{\wt{M}}}_t^{(\beta)}(\theta)dP_{\theta},
$$
which proves the validity of $(M_{\theta})$ for 
$A\in\mathcal{K}_s^{\sagp}$; hence step (d) follows.
	
\noindent{\bf (e)}
There exists a $P_{\vT}$-null set ${\qans}\in\mf{B}(\vY)$ such that for any 
$\theta\notin{\qans}$ condition $(M_{\theta})$ is satisfied, 
$M^{(\beta)}(\theta)$ is a 
martingale in $\mathcal{L}^1(P_{\theta})$,
and $\E_{P_{\theta}}[{\wt{M}}_t^{(\beta)}(\theta)]=1$.

To show (e), first note that the consistency of $\{Q_{\theta}\}_{\theta\in\vY}$ with $\vT$ yields the existence of a $Q_{\vT}$-null set 
${\qcns}\in\mf{B}(\vY)$ such that for any $\theta\notin{\qcns}$ and $B\in\mf{B}(\vY)$ we have $Q_{\theta}(\vT^{-1}(B))=\chi_{B}(\theta)$. 
In the same reasoning, there exists a $P_{\vT}$-null set ${\pcns}\in\mf{B}(\vY)$ such that for any $\theta\notin{\pcns}$ and $B\in\mf{B}(\vY)$ we have 
$P_{\theta}(\vT^{-1}(B))=\chi_{B}(\theta)$.
Then, putting ${\qans}:={\newd}\cup{V}\cup{\qcns}\in\mf{B}(\vY)$,
we get $P_{\vT}({\qans})=0$. 

Next denote by $\mcol{{\mathcal{K}}_s}$ the family of all $A\in{\haf}_s$ satisfying $(M_{\theta})$ for any fixed $\theta\notin{\qans}$. 
We shall first show that
\begin{equation}\label{incl}
\haf_s^{\sagp}\cup\sigma(\vT)\subseteq\mcol{{\mathcal{K}}_s}.
\end{equation}

To do it, let $A=\vT^{-1}(B)$ for $B\in\mf{B}(\vY)$. We then get for any 
fixed $\theta\in{(\qans)^{c}}\cap B$ that 
$$
\E_{P_{\theta}}[\chi_A{\wt{M}}_t^{(\beta)}(\theta)]
=\E_{P_{\theta}}[\chi_{\vT^{-1}(B)}{\wt{M}}_t^{(\beta)}(\theta)]=1=Q_{\theta}(\vT^{-1}(B)),
$$
where the second equality follows by the consistency of $\{P_{\theta}\}_{\theta\in\vY}$ with $\vT$.

Applying the same reasoning we get $\E_{P_{\theta}}[\chi_A{\wt{M}}_t^{(\beta)}(\theta)]=0=Q_{\theta}(\vT^{-1}(B))$ 
for any fixed $\theta\in{(\qans)^{c}}\cap B^c$; hence $A\in\mcol{{\mathcal{K}}_s}$ and so $\sigma(\vT)\subseteq{{\mathcal{K}}_s}\neq\emptyset$.

Note now that the family 
$\mathcal{D}_s$ of all $A\in\haf_s^{\sagp}$ satisfying $(M_{\theta})$ for any $\theta\notin{\qans}$ can be easily shown 
to be a Dynkin class containing $\mathcal{K}_s^{\sagp}$; 
hence applying the Dynkin Lemma 
we obtain that 
$\haf_s^{\sagp}=\sigma(\mathcal{K}_s^{\sagp})\subseteq\mathcal{D}_s\subseteq\haf_s^{\sagp}$, implying 
$\haf_s^{\sagp}=\mathcal{D}_s$. 
Consequently, $\haf_s^{\sagp}\subseteq\mcol{{\mathcal{K}}_s}$,
which proves (\ref{incl}).
\smallskip

Define now the class
$$
{\mathcal{G}}:=\Bigl\{\bigcap_{k=1}^{m}A_k:
A_k\in\haf_s^{\sagp}\cup\sigma(\vT),\;m\in\N\biggr\}.
$$
Then the following condition holds true:
\begin{equation}\label{dinc}
{\mathcal{G}}\subseteq\mcol{{\mathcal{K}}_s}.
\end{equation}
To show it, fix on arbitrary $G\in{\mathcal{G}}$. Then there exist a number $m\in\N$ and a finite sequence $\{A_k\}_{k\in\{1,\ldots,m\}}$ in
$\haf_s^{\sagp}\cup\sigma(\vT)$ 
such that $G=\bigcap_{k=1}^{m}A_k$. Setting
$$
I_{\vT}:=\{k\in\{1,\ldots,m\}:A_k\in\sigma(\vT)\}\;\;\mbox{and}\;\;
I_H:=\{k\in\{1,\ldots,m\}:A_k\in\haf_s^{\sagp}\setminus\sigma(\vT)\},
$$
we get $I_{\vT}\cup I_H=\{1,\ldots,m\}$ 
as well as
\begin{equation}
\bigcap_{k\in I_{\vT}}A_k\in\sigma(\vT)\quad\mbox{and}\quad
\bigcap_{k\in I_H}A_k\in\haf_s^{\sagp}\adc{.}	
\label{tm1}	
\end{equation}
It then follows by the first part of (\ref{tm1}) that there exists a set $D\in\mf{B}(\vY)$ such that $\vT^{-1}(D)=\bigcap_{k\in I_{\vT}}A_k$, implying together with the second part of (\ref{tm1}), condition (\ref{incl}) and  the consistency of $\{Q_{\theta}\}_{\theta\in\vY}$ and $\{P_{\theta}\}_{\theta\in\vY}$ with $\vT$ that for any
fixed $\theta\in{(\qans)^{c}}\cap D$ we get
\begin{eqnarray*}
Q_{\theta}(G)&=&Q_{\theta}\Bigl(\vT^{-1}(D)\cap\Bigl(\bigcap_{k\in I_H}A_k\Bigr)\Bigr)
=Q_{\theta}\Bigl(\bigcap_{k\in I_H}A_k\Bigr)
=\E_{P_{\theta}}[\chi_{\bigcap_{k\in I_H}A_k}{\wt{M}}_s^{(\beta)}(\theta)]\\
&=&\int_{\vT^{-1}(D)}\chi_{\bigcap_{k\in I_H}A_k}{\wt{M}}_s^{(\beta)}(\theta)dP_{\theta}
=\E_{P_{\theta}}[\chi_G{\wt{M}}_s^{(\beta)}(\theta)]=\E_{P_{\theta}}[\chi_G{\wt{M}}_t^{(\beta)}(\theta)],
\end{eqnarray*}
while for any fixed $\theta\in{(\qans)^{c}}\cap D^c$ we clearly 
get that
$$
\E_{P_{\theta}}[\chi_G{\wt{M}}_t^{(\beta)}(\theta)]=0
=Q_{\theta}\Bigl(\vT^{-1}(D)\cap\bigl(\bigcap_{k\in I_H}A_k\bigr)\Bigr)=Q_{\theta}(G).
$$
Consequently, $G\in\mcol{{\mathcal{K}}_s}$; hence ${\mathcal{G}}\subseteq\mcol{{\mathcal{K}}_s}$, which proves (\ref{dinc}).
\smallskip

It can be easily seen that $\mathcal{K}_s$ is a Dynkin class; hence taking into account condition (\ref{dinc}), 
we may apply the Dynkin Lemma 
to obtain condition $\haf_s=\mathcal{K}_s$.

For any $\theta\notin{\qans}$ condition ($M_{\theta}$) 
implies
$$
\int_A{\wt{M}}_s^{(\beta)}(\theta)dP_{\theta}=\int_A{\wt{M}}_t^{(\beta)}(\theta)dP_{\theta}\quad\mbox{for any} \quad A\in\haf_s;
$$
hence the family $M^{(\beta)}(\theta)$ is a 
martingale in $\mathcal{L}^1(P_{\theta})$.

Again by ($M_{\theta}$) we obtain for $A=\vO$ that condition $\E_{P_{\theta}}[{\wt{M}}_t^{(\beta)}(\theta)]=Q_{\theta}(\vO)=1$ holds for any $\theta\notin{\qans}$ true, completing 
the proof of (e) and of the implication $(ii)\Longrightarrow (iii)$.
\smallskip

Ad $(iii)\Longrightarrow (i)$:  Assuming assertion $(iii)$ we get that the measures $P_{\theta}$ and $Q_{\theta}$ are progressively equivalent for all $\theta\notin{\qans}$. Thus, taking into account that $\{P_{\theta}\}_{\theta\in\vY}$ and $\{Q_{\theta}\}_{\theta\in\vY}$ are r.c.p.s of $P$ over $P_{\vT}$ and $Q$ over $Q_{\vT}$, respectively, consistent with $\vT$ we obtain assertion $(i)$.
\smallskip

Ad $(iii)\Longrightarrow (iv)$: It follows by assertion $(iii)$ that there exists an essentially unique $\beta\in\mathcal{F}_P$ satisfying 
condition $(**)$. Thus, we are left to show that there exists a $\xif\in\mathcal{R}_+(\vY)$ which is a Radon-Nikod\'{y}m derivative of $Q_{\vT}$ with respect to $P_{\vT}$ satisfying 
condition $(M_{\xif})$, and that the family $M^{(\beta)}(\vT)$ is a 
martingale in $\mathcal{L}^1(P)$ satisfying condition $\mathbb{E}_P[M_t^{(\beta)}(\vT)]=1$ for all $t\in\mathbb{R}_+$.
\smallskip

In fact, 
first note that the assumption  
$Q_{\vT}\sim P_{\vT}$, implies the existence of a $P_{\vT}$-a.s. positive Radon-Nikod\'{y}m derivative $\xif$ of $Q_{\vT}$ with  
respect to {$P_{\vT}$}
and write $\wt{M}_t^{\beta}(\vT):=e^{S_t^{(\beta)}(\vT)-t\vT(e^{\alpha(\vT)}-1)}$. Then condition 
\begin{equation}\label{ee1}
\E[\wt{M}_t^{(\beta)}(\vT)\mid\vT]=1\quad P\uph{\sigma(\vT)}-\mbox{a.s.}
\end{equation}
is valid, since 
\begin{eqnarray*}
\E[\wt{M}_t^{(\beta)}(\vT)\mid\vT]
&=&e^{-t\vT(e^{\alpha(\vT)}-1)}\sum_{n=0}^{\infty}e^{n\alpha(\vT)}P(\{N_t=n\}\mid\vT)\E_{P}[e^{\sum_{k=1}^{n}\gamma(X_k)}]\\
&=&e^{-t\vT(e^{\alpha(\vT)}-1)}e^{-t\vT}
\sum_{n=0}^{\infty}\frac{[(t\vT)e^{\alpha(\vT)}]^n}{n!}=1,
\end{eqnarray*}
where all equalities hold $P\uph\sigma(\vT)$-a.s. true.

Fix on arbitrary $A\in\haf_s$. Applying condition (\ref{ee1}) we obtain 
$$
\int\chi_AM_t^{(\beta)}(\vT)dP 
=\int\chi_A{\xif}(\vT)\wt{M}_t^{(\beta)}(\vT)dP
\leq\int{\xif}(\vT)\E_P[\wt{M}_t^{(\beta)}(\vT)\mid\vT]dP
=1.
$$
Thus, we may apply \cite{lm1}, Proposition 3.8 (ii),  
for $f=\vT$ and $u=(\chi_A\otimes\xif){\wt{M}}_t^{(\beta)}$, where $(\chi_A\otimes\xif)(\omega,\theta):=\chi_A(\omega)\xif(\theta)$ 
and $\wt{M}^{(\beta)}_t(\omega,\theta):=\wt{M}^{(\beta)}_t(\theta)(\omega)$ for all $(\omega,\theta)\in\vO\times\vY$, to get
$$
\int\chi_AM_t^{(\beta)}(\vT)dP=\int\chi_A{\xif}(\theta)\E_{P_{\theta}}[{\wt{M}}_t^{(\beta)}(\theta)]P_{\vT}(d\theta).
$$

The latter together with condition $(M_{\theta})$ yields
$$
Q(A)
=\int{Q}_{\theta}(A)Q_{\vT}(d\theta)
=\int\E_{P_{\theta}}[\chi_A{\wt{M}}_t^{(\beta)}(\theta)]\xif(\theta)P_{\vT}(d\theta)\\
=\int_AM_t^{(\beta)}(\vT)dP;
$$
hence
$M^{(\beta)}(\vT)$ is a 
martingale satisfying condition 
$\mathbb{E}_P[M_t^{(\beta)}(\vT)]=1$ for all $t\in\mathbb{R}_+$.
The implication $(iv)\Longrightarrow(i)$ is clear.
This completes the whole proof.\hfill$\Box$
\medskip

For given $P\in{\mathcal{M}}_{\sagp}$ and $Q\in\mathcal{M}_{\sagp,g}$ it follows that the measures $P$ and $Q$ are equivalent on each 
$\sigma$-algebra ${\haf}_t$. 
But this result does not in general hold true on ${\haf}_{\infty}$ as the following proposition shows. 

\begin{prop}\label{46}
Let be given $Q\in\mathcal{M}_{\sagp,g}$ and a r.c.p. $\{Q_{\theta}\}_{\theta\in\vY}$ of $Q$ over $Q_{\vT}$ consistent with $\vT$. Assume that $P_{\theta}\neq Q_{\theta}$ for $P_{\vT}$-a.a. $\theta\in\vY$. Then the measures $P$ and $Q$ are singular on 
${\haf}_{\infty}$, i.e. there exists a set $E\in\haf_{\infty}$ such that $P(E)=0$ if and only if $Q(E)=1$. 
\end{prop}

{\bf Proof.} 
Proposition \ref{12} together with Proposition \ref{tsp}, $(i)\Longrightarrow(ii)$ implies
that there exists a $P_{\vT}$- and $Q_{\vT}$-null set $\wt{H}\in\mf{B}(\vY)$ such that $S$ is a $P_{\theta}$-CPP($\theta,P_{X_1}$) 
and a $Q_{\theta}$-CPP($g(\theta),Q_{X_1}$)
for any $\theta\notin\wt{H}$.

By assumption there exists a $P_{\vT}$-null set 
$\wh{H}\in\mf{B}(\vY)$ such that $P_{\theta}\neq Q_{\theta}$ for any $\theta\notin\wh{H}$. Note that we may and do assume that 
$\wh{H}$ contains $\wt{H}$.

Let us fix on arbitrary $\theta\notin\wh{H}$. Since $P_{\theta}\neq Q_{\theta}$, it follows by 
\cite{dh}, Proposition 2.2,
that the measures 
$P_{\theta}$ and $Q_{\theta}$ are singular on ${\haf}_{\infty}^S$, 
implying that there exists a set $E\in{\haf}_{\infty}$ such that 
$$
P_{\theta}(E)=0\Longleftrightarrow Q_{\theta}(E)=1.
$$
But the latter together with Proposition \ref{tsp}, $(i)\Longrightarrow(ii)$, yields that
$$
P(E)=\int{P}_{\theta}(E)P_{\vT}(d\theta)=0\Longleftrightarrow Q(E)=\int{Q}_{\theta}(E)Q_{\vT}(d\theta)=1.
$$
Consequently, the measures $P$ and $Q$ are singular on ${\haf}_{\infty}$.\hfill$\Box$ 

\section{The Characterization}\label{lm35}

In the next proposition a construction of non-trivial probability spaces admitting CMPPs is given. Such a construction not only 
it is of vital importance for the main result of this work (i.e. Theorem \ref{qQ}) but also extends a similar construction for 
MPPs, see \cite{lm5jmaslo}, Theorem 3.1. 

To prove it, we recall the following notations concerning products of probability spaces.
By $(\vO\times\vY,\vS\otimes{T},P\otimes{Q})$ is denoted the product probability space of $(\vO,\vS,P)$ and $(\vY,T,Q)$, 
and by $\pi_{\vO}$ and $\pi_{\vY}$ the 
canonical projections from $\vO\times\vY$ onto $\vO$ and $\vY$, respectively. 
If $f$ is a real-valued function defined on $\vO\times\vY$, then we shall be using the ordinary notation $f_{\omega}$ and $f^{y}$ for the functions obtained from $f$ by fixing $\omega$ and $y$, respectively. In a similar way the sets $C_{\omega}$, $C^{y}$ being sections of a set $C\subseteq\vO\times\vY$ are defined. 
If $(\vO,\vS,P)$ is a probability space and $I$ a non-empty index set, we write $P_I$ for the product measure on $\vO^I$ and $\vS_I$ for its domain.

{\em Throughout what follows, we put $\wt\vO:=\vY^{\N}\times\vY^{\N}$, 
$\wt\vS:=\mf{B}(\wt\vO)=\mf{B}(\vY)_{\N}$, $\vO:=\wt\vO\times\vY$ and $\vS:=\mf{B}(\vO)$ for simplicity}.
\smallskip
  
\begin{prop}\label{conb}
Let be given two probability measures $\mu$ and $\rho$ on $\mf{B}(\vY)$ and a positive $\mf{B}(\vY)$-measurable function ${h}$.
Then there exist:
\begin{enumerate}
\item 
A family $\{P_{\theta}\}_{\theta\in\vY}$ of probability measures 
$P_{\theta}:=(\mathbf{Exp}({h}(\theta)))_{\N}\otimes\mu_{\N}\otimes\delta_{\theta}$, 
where $\delta_{\theta}$ is the Dirac measure on $\mf{B}(\vY)$ concentrated on $\theta$, and a probability measure $P$ on $\vS$ 
such that $\{P_{\theta}\}_{\theta\in\vY}$ is a r.c.p. of 
$P$ over $\rho$ consistent with $\vT:=\pi_{\vY}$, where $\pi_{\vY}$ is the canonical projection from $\vO$ onto $\vT$,
and $P_{\vT}=\rho$.
\item
A counting process $\cnp$ and a \cspw $\csp$ such that $P_{X_1}=\mu$ and the quadruplet
$(P,\cip,\csp,\vT)$, where $\cip$ is the interarrival process associated with $\cnp$, satisfies conditions \textnormal{(a1)}
and \textnormal{(a2)}, and the pair $\rp$ is both a $P$- and a $P_{\theta}$-risk process for any $\theta\in\vY$,
inducing an \agpw $\agp$ being a $P$-CMPP$({h}(\vT),\mu)$.
\end{enumerate}
\end{prop}

{\bf Proof.} For any $\theta\in\vY$ put $\nu_{\theta}:=\mathbf{Exp}({h}(\theta))$ and $\wt{P}_{\theta}:=(\nu_{\theta})_{\N}\otimes\mu_{\N}$. 

\noindent{\bf(a)} 
Each set-function $\wt{P}_{\theta}$ is a probability measure on $\wt\vS$ and for any fixed $\wt{E}\in\wt{\vS}$ the function 
$\theta\longmapsto\wt{P}_{\theta}(\wt{E})$ is $\mf{B}(\vY)$-measurable. 

In fact, it is clear that $\wt{P}_{\theta}$ is for any fixed $\theta\in\vY$ a probability measure on $\wt{\vS}$. 
Furthermore, since $\nu_{\theta}(D)=\int\chi_D(x){h}(\theta)e^{-x{h}(\theta)}\leb(dx)$ for each $D\in\mf{B}(\vY)$, where $\leb$ denotes the restriction of the Lebesgue measure to $\mf{B}(\vY)$,
it follows by a monotone class argument that the function 
$\theta\longmapsto\nu_{\theta}(D)$ is $\mf{B}(\vY)$-measurable for any fixed $D\in\mf{B}(\vY)$; hence
${\wt{P}}_{\bdot}(A\times B)$ is a $\mf{B}(\vY)$-measurable function 
for any fixed $A\times B\in\mf{B}(\vY^{\N})\times\mf{B}(\vY^{\N})$.
It then follows again by a monotone class argument \mcol{that} the function $\theta\longmapsto\wt{P}_{\theta}(\wt{E})$ is $\mf{B}(\vY)$-measurable 
for any fixed $\wt{E}\in\wt\vS$.

\noindent{\bf(b)} Define the set-functions 
$\wt{P}:\wt\vS\longrightarrow[0,\infty]$ and $Q:\vS\longrightarrow[0,\infty]$ 
by means of
$$
\wt{P}(\wt{E}):=\int\wt{P}_{\theta}(\wt{E}){\rho}(d\theta)
\quad\text{for all}\quad \wt{E}\in\wt\vS,
$$
and
$$
P(E):=\int\wt{P}_{\theta}(E^{\theta}){\rho}(d\theta)\quad\text{for all}\quad{E}\in\vS.
$$
Then $\wt{P}$ and $P$ are probability measures on $\wt\vS$ and $\vS$, respectively, such that 
$\{P_{\theta}\}_{\theta\in\vY}$ is a r.c.p. of $P$ over {$\rho$} consistent with $\vT$ and $P_{\vT}=\rho$.

In fact, obviously $\wt{P}$ and $P$ are probability measures. 
Furthermore, for each $B\in\mf{B}(\vY)$ we have $P_{\vT}(B)=\int\wt{P}_{\theta}([\wt\vO\times{B}]^{\theta})\rho(d\theta)
=\int_B\wt{P}_{\theta}(\wt\vO)\rho(d\theta)=\rho(B)$. 
Since $\wt{P}$ and $P$ are the marginals of $P$ and by step (a) the function $\theta\longmapsto\wt{P}_{\theta}(\wt{E})$ is $\mf{B}(\vY)$-measurable 
for any fixed $\wt{E}\in\wt\vS$, we obtain that $\{\wt{P}_{\theta}\}_{\theta\in\vY}$ is a product r.c.p. on $\wt\vS$ for $P$ with respect to {$\rho$} 
(see \cite{mms3}, page 2390 for the definition and its properties). 
Put $P_{\theta}:=\wt{P}_{\theta}\otimes\delta_{\theta}$, where $\delta_{\theta}$ is the Dirac measure on $\vY$.
It then follows by \cite{lm5jmaslo}, Lemma 2.4, that 
$\{P_{\theta}\}_{\theta\in\vY}$ is a r.c.p. of $P$ over 
{$\rho$} consistent with $\vT$. This completes the proof of $(i)$. 

\noindent{\bf(c)}
There exists a \cnpw $\cnp$ and a 
\cspw $\csp$ 
such that the quadruplet $(P,\cip,\csp,\vT)$ satisfies conditions \textnormal{(a1)}
and \textnormal{(a2)}, and the pair $\rp$ is both a $P$- and a $P_{\theta}$-risk process for any $\theta\in\vY$.

In fact, first fix on arbitrary $n\in\N$ and $\theta\in\vY$. Denote by $\pi_{\wt\vO}$ 
the canonical projection from $\vO$ onto $\wt\vO$ and by $\wt{W}_n$ and $\wt{X}_n$ the canonical projections from $\wt\vO$ onto the
$n$-coordinate of the first and the second factor of $\vO=\vY^{\N}\times\vY^{\N}\times\vY$, respectively. 
Putting $W_n=\wt{W}_n\circ\pi_{\wt\vO}$ and $X_n=\wt{X}_n\circ\pi_{\wt\vO}$ we get 
\begin{equation}\label{mast2} 
\nu_{\theta}={{P}}_{\theta}\circ{W_n}^{-1}=(\wt{P}_{\theta})_{\wt{W}_n}\quad\mbox{and}\quad 
(P_{\theta})_{X_n}=({\wt{P}}_{\theta})_{\wt{X}_n}=\mu.
\end{equation}
Since $(\wt\vO,\wt\vS,\wt{P}_{\theta})$ is a product probability space for any $\theta\in\vY$, and since $\wt{W}_n, \wt{X}_n$ are 
the canonical projections, it follows by standard arguments that the processes $\{\wt{W}_n\}_{n\in\N}$ and $\{\wt{X}_n\}_{n\in\N}$ are $\wt{P}_{\theta}$-independent and $\wt{P}_{\theta}$-mutually independent; hence the processes $W$ and $X$ are $P_{\theta}$-independent and $P_{\theta}$-mutually independent for any $\theta\in\vY$.

Let $\cnp$ be the counting process associated with $\cip$ (cf. e.g. \cite{sch}, page 6 and Theorem 2.1.1). 
The fact that $W$ and $X$ are $P_{\theta}$-mutually independent together with (\ref{Z0}) implies that
the processes $\cnp$ and $\csp$ are mutually independent under $P_{\theta}$.
Thus the pair $\rp$ is a $P_{\theta}$-risk process. 

Since $\{{P}_{\theta}\}_{\theta\in\vY}$ is a r.c.p. of $P$ over $\rho$ consistent with $\vT$ by (b), 
applying \cite{lm1}, Lemma 4.1, 
we get that $\cip$ and $\csp$ are $P$-conditionally independent, 
i.e. that $P$ satisfies condition \textnormal{(a1)}.
Furthermore, 
condition (\ref{mast2}) together with the fact that $\{{P}_{\theta}\}_{\theta\in\vY}$ is a r.c.p. of $P$ over 
$\rho$ consistent with $\vT$ by (b), yields $P_{X_n}=\mu$; hence condition \textnormal{(a2)} is satisfied by $P$.
Conditions \textnormal{(a1)} and \textnormal{(a2)} for $P$ together with the fact that $\{{P}_{\theta}\}_{\theta\in\vY}$ is a r.c.p. of $P$ over 
$\rho$ consistent with $\vT$, imply the 
$P$-mutual independence of $\cip$ and $\csp$.

The latter together with (\ref{Z0}) implies that the processes $\cnp$ and $\csp$ are
$P$-mutually independent.
Since $\csp$ is $P_{\theta}$-i.i.d. and $P$ satisfies condition \textnormal{(a2)}, 
then according to Lemma \ref{11}, $(ii)$ the process $\csp$ is $P$-i.i.d..
Thus the pair $\rp$ is a $P$-risk process. This completes the proof of (c).

\noindent{\bf(d)}
The process $\agp$ induced by the $P$-risk process $\rp$ is a $P$-CMPP$({h}(\vT),\mu)$.
In fact, since for all $\theta\in\vY$ the sequence $W$ is $P_{\theta}$-independent and $(P_{\theta})_{W_n}=\mathbf{Exp}(h(\theta))$ for all $n\in\N$  by (c), it follows that
$\cnp$ is a $P_{\theta}$-PP(${h}(\theta)$) 
(cf. e.g. \cite{sch}, Theorem 2.3.4). The latter together with the fact that, due to (c), the pair $\rp$ is a $P_{\theta}$-risk 
process implies that the induced \agpw $\agp$ is a $P_{\theta}$-CPP(${h}(\theta),\mu$) for any $\theta\in\vY$, 
implying according to Proposition \ref{12} that $\agp$ is a $P$-CMPP$({h}(\vT),\mu)$. Thus, assertion $(ii)$ follows, completing 
in this way the whole proof.\hfill$\Box$
\smallskip

Note that since $\vO$ is a Polish space, it follows that for each probability measure 
$P$ on ${\haf}_{\infty}$ 
there always exists a r.c.p. $\{P_{\theta}\}_{\theta\in\vY}$ of $P$ over $P_{\vT}$ consistent with $\vT$ (see Remark \ref{magd}(a)). 
\smallskip

{\em From now on, unless it is stated otherwise, we consider the probability space $(\vO,\vS,P)$ together with 
${h}:=id_{\vY}$, the random variable $\vT$,
the r.c.p. $\{P_{\theta}\}_{\theta\in\vY}$ as well as the processes $\cnp$, $\cip$, $\csp$ and $\agp$ constructed in Proposition \ref{conb}}.
 
\begin{rem}\label{n49}
\normalfont
Such an assumption concerning $(\vO,\vS,P)$ is not a restrictive one,
since our interest does not exceed the information generated by the \agpw and the structure parameter. In fact, since 
$\vS=\sigma(\{\cip,\csp,\vT\})$, we clearly get that $\haf_{\infty}\subseteq\vS$. But due to 
\cite{sch}, Lemma 2.1.2 and Remark \ref{qS1}, (a), we have something more, i.e. that $\vS={\haf}_{\infty}$. More precisely, 
Lemma 2.1.2 of \cite{sch} yields 
$\sigma(\cip)\subseteq\haf_{\infty}^{\sagp}$, 
while by Remark \ref{qS1}, (a), we get 
$X_n^{-1}(B)\cap\{T_n\leq l\}\in\haf_l^{\sagp}$ for all $n,l\in\N_0$ and $B\in\mf{B}(\vY)$, implying $X_n^{-1}(B)\in\haf_{l}^{\sagp}$, 
since 
$\bigcup_{l\in\N}\{T_n\leq l\}=\vO$ 
by our assumption that there cannot be infinite claims in finite time; hence 
$\sigma(\csp)\subseteq\haf_{\infty}^{\sagp}$ 
and so $\vS\subseteq{\haf}_{\infty}$ follows.
\end{rem}

\begin{thm}\label{qQ} 
The following holds true:
\begin{enumerate}
\item
For each pair $(g,Q)\in\mf{M}_+(\vY)\times\mathcal{M}_{\sagp,g}$ 
there exists an essentially unique pair 
$(\beta,\xif)\in\mathcal{F}_P\times\mathcal{R}_+(\vY)$ with $\xif$ 
being a Radon-Nikod\'{y}m derivative of $Q_{\vT}$ with respect to $P_{\vT}$, satisfying conditions $(**)$ and  $(M_{\xif})$, and such that the family $M^{(\beta)}(\vT)$, involved in $(M_{\xif})$, is a 
martingale in $\mathcal{L}^1(P)$; 
\item 
Conversely, for every pair $(\beta,\xif)\in\mathcal{F}_P\times\mathcal{R}_+(\vY)$ there exists a unique pair 
$(g,Q)\in\mf{M}_+(\vY)\times\mathcal{M}_{\sagp,g}$  satisfying conditions $(**)$ and  $(M_{\xif})$, and such that the family $M^{(\beta)}(\vT)$, involved in $(M_{\xif})$, is a 
martingale in $\mathcal{L}^1(P)$; 
\item
In both cases $(i)$ and $(ii)$ there exist an essential unique 
r.c.p. $\{Q_{\theta}\}_{\theta\in\vY}$ of $Q$ over $Q_{\vT}$ consistent with $\vT$, and a $P_{\vT}$-null set 
${\qans}\in\mf{B}(\vY)$ such that for any $\theta\notin{\qans}$ 
conditions $Q_{\theta}\in\wt{\mathcal{M}}_{\sagp,g}$ and $(M_{\theta})$ are valid, and such that the family $\wt{M}^{(\beta)}(\theta)$, involved in 
$(M_{\theta})$, is a 
martingale in $\mathcal{L}^1(P_{\theta})$. 
\end{enumerate}
\end{thm}

\noindent {\bf Proof.}
Assume that $Q\in\mathcal{M}_{\sagp,g}$. There always exists an essentially unique r.c.p. $\{Q_{\theta}\}_{\theta\in\vY}$ of $Q$ over $Q_{\vT}$ consistent with $\vT$ (see Remark \ref{magd}(a)). So, we may apply Proposition \ref{tsp} to derive assertion $(i)$.

\noindent Ad $(ii)$: Let $(\beta,\xif)\in\mathcal{F}_P\times\mathcal{R}_+(\vY)$. 
Define the function $g:\vY\longrightarrow\vY$ by means of $g(\theta):=\theta{e}^{\alpha(\theta)}$ for each $\theta\in\vY$. 

Fix on arbitrary $t\in\mathbb{R}_+$ and $\theta\in{\vY}$, 
and define the set-functions $\wt\nu_{\theta}, \wt\mu, \wt{R}:\mf{B}\longrightarrow [0,\infty]$ by means of 
$$
\wt\nu_{\theta}:=\wt\nu_{\theta, g}:=\mathbf{Exp}(g(\theta)),
$$
$$
\wt\mu(B):=\int_B{e}^{\gamma(x)}P_{X_1}(dx)\quad\mbox{for any}\quad B\in\mf{B}(\vY),
$$
$$
\wt{R}(E):=\E_P[\chi_{\vT^{-1}(E)}\xif(\vT)]\quad\mbox{for any}\quad E\in \mf{B}(\vY),
$$
respectively. Clearly, $\wt{R}$ and $\wt\mu$ are are probability measures on $\mf{B}(\vY)$, since by assumption $\E_P[\xif(\vT)]=1$ and $\E_P[e^{\gamma(X_1)}]=1$, respectively. 
Thus, applying Proposition \ref{conb} for $g, \wt\nu_{\theta}, \wt\mu$ and $\wt{R}$ in the place of 
${h},\nu_{\theta}, \mu$ and $\rho$, respectively, we obtain a family $\{Q_{\theta}\}_{\theta\in\vY}$ of probability measures 
$Q_{\theta}:=\mathbf{Exp}(g(\theta))\otimes\wt{R}_{\N}\otimes\delta_{\theta}$ on $\vS$ and a unique probability measure $Q$ on $\vS$ satisfying assumptions 
\textnormal{(a1)} and \textnormal{(a2)}, and such that $\{Q_{\theta}\}_{\theta\in\vY}$ is a r.c.p. of $Q$ over $Q_{\vT}=\wt{R}$ consistent with $\vT$ and $S$ is a 
$Q$-CMPP$(g(\vT),Q_{X_1})$ with $Q_{X_1}=\wt\mu$. 
The latter implies that $Q_{\vT}\sim P_{\vT}$ and $Q_{X_1}\sim P_{X_1}$. Thus, according to Proposition \ref{tsp}, we get $Q\stackrel{pr}{\sim}P$,  
implying both $Q\in\mathcal{M}_{\sagp,g}$ and the conclusion of $(ii)$. 

\noindent Ad $(iii)$: 
Since $\vO$ is a Polish space, in both cases there exists an essentially unique r.c.p. $\{Q_{\theta}\}_{\theta\in\vT}$ of $Q$ over $Q_{\vT}$ consistent with $\vT$ (see again Remark \ref{magd}(a)); hence assertion $(iii)$ is an immediate consequence of Proposition \ref{tsp}. This completes the proof.\hfill$\Box$

\begin{rems}\label{qr}
\normalfont
{\bf (a)} 
Proposition 2.2 of Delbaen \& Haezendonck \cite{dh} can be derived as a special case of Theorem \ref{qQ} if the probability distribution of the random variable $\vT$ is degenerate \mcol{under $P$} at some $\theta_0>0$.

{\bf (b)} 
Comparing our main result (Theorem \ref{qQ}) with the corresponding
characterization of Meister \cite{mei}, Corollary 2.10, some major differences that should be pointed out are spotted:
Our result is proven without either the completeness assumption of the underlying filtration (see \cite{mei}, Assumption 1.4), or the square integrability of the random variables
$\vT$ and $X_n$ ($n\in\N$) (see \cite{mei}, Section 1.5), or the $P$-independence assumption of the 
sequences $\cip$ and $\csp$ (see \cite{mei}, Section 2.2).
The latter together with the fact that
in Theorem \ref{qQ} we consider a simpler and more natural underlined probability space, constructed in Proposition \ref{conb},    
also containing the information of the structural parameter $\vT$, allows us to obtain by means of Theorem \ref{qQ} 
not only a proper measure $Q$ as happens in \cite{mei}, Corollary 2.10, but also
a whole family of proper measures $Q_{\theta}$ along with $Q$. 
\end{rems}

\begin{ntr}\label{nr}
\normalfont
{\bf (a)} 
For given $g\in\mf{M}_+(\vY)$ define the classes 
$\mathcal{M}^{\ell}_{\sagp,g}:=\{Q\in\mathcal{M}_{\sagp,g}: \E_Q[X^{\ell}_1]<\infty\}$ and 
$\mathcal{F}_P^{\ell}:=\{\beta\in\mathcal{F}_P: \E_P[X^{\ell}_1e^{\gamma(X_1)}]<\infty\}$ for $\ell\in\{1,2\}$. 
For $g=id_{\vY}$ put $\mathcal{M}^{\ell}_{\sagp}:=\mathcal{M}^{\ell}_{\sagp,g}$. It can be easily seen that Theorem \ref{qQ} remains true, if we replace the classes   
$\mathcal{M}_{\sagp,g}$ and $\mathcal{F}_P$ by their subclasses $\mathcal{M}^{\ell}_{\sagp,g}$ and $\mathcal{F}_P^{\ell}$ for $\ell=1,2$, respectively.

{\bf (b)} 
For any $\theta\in\vY$ denote by $\wt{\mathcal{M}}^{\ell}_{\sagp,g}:=\wt{\mathcal{M}}^{\ell}_{\sagp,g}(\theta)$ the class of all probability measures 
$Q_\theta\in\wt{\mathcal{M}}_{\sagp,g}$ such that $\E_{Q_{\theta}}[X_1^{\ell}]<\infty$. 

{\bf (c)}
The probability measure $P$ constructed in Proposition \ref{conb} is an element of $\mathcal{M}_S$, 
and assuming that $\int_{\vY}x^{\ell}\mu(dx)<\infty$ for $\ell\in\{1,2\}$ we get $\E_P[X_1^{\ell}]<\infty$; hence $P\in\mathcal{M}_S^{\ell}$.

{\bf (d)}
For given $g\in\mf{M}_+(\vY)$ and $\ell\in\{1,2\}$ define the class $\mathcal{R}^{*,\ell}_+(\vY)$ 
of all $\xif\in{\mathcal{R}}_+(\vY)$ such that $\E_P[\xif(\vT)g^{\ell}(\vT)]<\infty$ as well as
the class of $\mathcal{M}_{\sagp,g}^{*,\ell}$ of all measures $Q\in\mathcal{M}^{\ell}_{\sagp,g}$ 
such that $\E_Q[g^{\ell}(\vT)]<\infty$. We then get that 
$\mathcal{M}^{*,\ell}_{\sagp,g}\subseteq\mathcal{M}^{\ell}_{\sagp,g}\subseteq\mathcal{M}_{\sagp,g}$. 
For $g=id_{\vY}$ put $\mathcal{M}^{*,\ell}_{\sagp}:=\mathcal{M}_{\sagp,g}^{*,\ell}$. 
\end{ntr}

In the next result we give a characterization of all elements of $\mathfrak{M}_+(\vY)\times\mathcal{M}^{*,\ell}_{\sagp,g}$ in terms of the elements of
$\mathcal{F}_P^{\ell}\times\mathcal{R}^{*,\ell}_+(\vY)$. 

\begin{cor}\label{qm0}
Let $\ell\in\{1,2\}$ 
and $P\in\mathcal{M}_{\sagp}^{*,\ell}$. The following holds true:

\begin{enumerate}
\item %
For each pair $(g,Q)\in\mf{M}_+(\vY)\times\mathcal{M}^{*,\ell}_{\sagp,g}$ 
there exists an essentially unique pair $(\beta,\xif)\in\mathcal{F}_P^{\ell}\times\mathcal{R}^{*,\ell}_+(\vY)$ with $\xif$ 
being a Radon-Nikod\'{y}m derivative of $Q_{\vT}$ with respect to $P_{\vT}$, satisfying conditions $(**)$ and $(M_{\xif})$; 
\item 
Conversely, for each pair $(\beta,\xif)\in\mathcal{F}_P^{\ell}\times\mathcal{R}^{*,\ell}_+(\vY)$ there exists 
a unique pair $(g,Q)\in\mf{M}_+(\vY)\times\mathcal{M}^{*,\ell}_{\sagp,g}$ satisfying conditions $(**)$ and $(M_{\xif})$; 
\item
In both cases there exist an essentially unique r.c.p. $\{Q_{\theta}\}_{\theta\in\vY}$ of $Q$ over $Q_{\vT}$ consistent with $\vT$ and a $P_{\vT}$-null set 
${\qans}\in\mf{B}(\vY)$ such that for any $\theta\notin{\qans}$ conditions $Q_{\theta}\in\wt{\mathcal{M}}^{\ell}_{\sagp,g}$ and $(M_{\theta})$ are fulfilled. 
\end{enumerate}
\end{cor}
{\bf Proof.} 
Ad $(i)$: Assume that $(g,Q)\in\mf{M}_+(\vY)\times\mathcal{M}_{\sagp,g}^{*,\ell}$. 
Since clearly $\mathcal{M}_{\sagp,g}^{*,\ell}\subseteq\mathcal{M}_{\sagp,g}$, we get by Theorem \ref{qQ}, $(i)$ that there exist an essentially unigue pair $(\beta,\xif)\in\mathcal{F}_P\times\mathcal{R}_+(\vY)$ with $\xif$ being a Radon-Nikod\'{y}m derivative of $Q_{\vT}$ with respect to $P_{\vT}$ satisfying conditions $(*)$ and $(M_{\xif})$. 
But since $Q\in\mathcal{M}_{\sagp,g}^{*,\ell}$ we have that 
$g^{\ell}(\vT)\in\mathcal{L}^1(Q)$, implying that 
$(\beta,\xif)\in\mathcal{F}_P^{\ell}\times\mathcal{R}_{+}^{*,\ell}(\vY)$. 

Ad $(ii)$: Let 
$(\beta,\xif)\in\mathcal{F}_P^{\ell}\times\mathcal{R}^{*,\ell}_+(\vY)$.
Since $\mathcal{F}_P^{\ell}\subseteq\mathcal{F}_P$ it follows by Theorem \ref{qQ}, $(ii)$, that there exists a unique pair $(g,Q)\in\mf{M}_+(\vY)\times\mathcal{M}_{\sagp,g}$ satisfying conditions $(*)$ and $(M_{\xif})$ 
Furthermore, assumption $(\beta,\xif)\in\mathcal{F}_P^{\ell}\times\mathcal{R}^{*,\ell}_+(\vY)$ implies both $X_1^{\ell}\in\mathcal{L}^1(Q)$ and $g^{\ell}(\vT)\in\mathcal{L}^1(Q)$; hence $(g,Q)\in\mf{M}_+(\vY)\times\mathcal{M}^{*,\ell}_{\sagp,g}$. 

Ad $(iii)$: 
Due to Theorem \ref{qQ} there exist an essentially unique r.c.p. $\{Q_{\theta}\}_{\theta\in\vY}$ 
of $Q$ over $Q_{\vT}$ consistent with $\vT$ and a $P_{\vT}$-null set ${\qans}\in\mf{B}(\vY)$ such that for any 
$\theta\notin{\qans}$ condition $(M_{\theta})$ is fulfilled. Condition $Q_{\theta}\in\wt{\mathcal{M}}^{\ell}_{\sagp,g}$ for any $\theta\notin{\qans}$ follows as in the proof of $(i)$. 
This completes the proof.\hfill$\Box$
\medskip 

\section{Compound mixed Poisson processes and martingales}\label{cmppm}

In this section applying our results we find out a 
wide class of {\em canonical} processes satisfying the condition 
of {\em no free lunch with vanishing risk} 
(written (NFLVR) for short) (see \cite{ds}, Definition 8.1.2). 
Let $\ell\in\{1,2\}$. For given real-valued process 
$\Zp:=\{Z_t\}_{t\in{\R}_+}$ 
a process on $(\vO,\vS)$ by 
$\M^{\ell}_{S,g}(\Zp):=\M^{\ell}_{P,S,g}(\Zp)$ 
will be denoted the class of all probability measures $Q\in\mathcal{M}_{\sagp,g}^{\ell}$ such that $\Zp$ is a 
martingale in $\mathcal{L}^{\ell}(Q)$.
The elements of $\M^{\ell}_{S,g}(\Zp)$ 
are called {\bf progressively equivalent martingale measures} (with respect to $\Zp$). 

We say that the process $\Zp$ satisfies condition (PEMM) if 
$\M^{\ell}_{S,g}(\Zp)\not=\emptyset$. 
Moreover, let $T>0, \T:=[0,T], R_T:=R\uph\haf_T$, where $R$ is any measure on $\vS$, $\Zp_{\T}:=\{Z_t\}_{t\in\T}$ 
and $\haf_{\T}:=\{\haf_t\}_{t\in\T}$. 

By $\M^{\ell}_e(\Zp_{\T}):=\M^{\ell}_{e,P,S,g}(\Zp_{\T})$ 
will be denoted the class of all probability measures $Q_T$ on $\haf_T$ such that 
$Q_T\sim{P}_T$, $S_{\T}$ is a $Q_T$-CMPP($g(\vT),(Q_T)_{X_1})$ and $\Zp_{\T}$ is a $(Q_T,\haf_{\T})$-martingale in 
$\mathcal{L}^{\ell}(Q_T)$.
The elements of $\M_e^{\ell}(\Zp_{\T})$ are called {\bf equivalent martingale measures} (with respect to $\Zp_{\T}$).
We say that the process $\Zp_{\T}$ satisfies condition (EMM) if $\M_e^{\ell}(\Zp_{\T})\not=\emptyset$. 

\begin{prop}\label{cqr}
Let $\ell\in\{1,2\}$, let $g\in\mathfrak{M}_+(\vY)$ and let $P\in{\mathcal{M}}^{\ell}_{\sagp}$. 
\begin{enumerate}
\item 
If $Q\in\mathcal{M}^{\ell}_{\sagp,g}$ then 
there exists a $Q_{\vT}$-null set ${\moic}\in\mf{B}(\vY)$ such that for any $\theta\notin{\moic}$ the process 
$\Vp(\theta):=\{\Vp_t(\theta)\}_{t\in\R_+}:=\{S_t-\E_{Q_{\theta}}[S_t]\}_{t\in\R_+}$ 
is a 
martingale in $\mathcal{L}^{\ell}(Q_{\theta})$; 
\item
$\mathcal{M}_{\sagp,g}^{*,\ell}=\M^{\ell}_{S,g}(\Vp(\vT))$, 
where $\Vp(\vT):=\{\Vp_t(\vT)\}_{t\in\R_+}:=\{S_t-tg(\vT)\E_P[X_1e^{\gamma(X_1)}]\}_{t\in\R_+}$;
\item
if $Q\in\mathcal{M}^{*,\ell}_{\sagp,g}$ then the process 
$\Vp:={\{\Vp_t\}_{t\in\R_+}:=}\{S_t-\E_{Q}[S_t]\}_{t\in\R_+}$ is a 
martingale in $\mathcal{L}^{\ell}(Q)$ if and only if the probability distribution of 
$g(\vT)$ is degenerate at $g(\theta_0)$ for some $\theta_0>0$.
\end{enumerate}
\end{prop}

{\bf Proof.} Let us fix on arbitrary 
$\ell\in\{1,2\}$ and let $Q\in\mathcal{M}_{\sagp,g}$. Then 
according to Proposition \ref{12}  there exists a $Q_{\vT}$-null set ${\lcnsq}\in\mf{B}(\vY)$ such that for any $\theta\notin{\lcnsq}$ the family $\agp$ is a $Q_{\theta}$-CPP$(g(\theta),(Q_{\theta})_{X_1})$, which yields that $\cnp$ is a 
 $Q_{\theta}$-PP($g(\theta)$), 
implying together with \cite{sch}, Theorem 5.1.3, that the family $\agp$ has $Q_{\theta}$-stationary independent increments. 

Ad $(i)$: If $Q\in\mathcal{M}^{\ell}_{\sagp,g}$ then the latter together with the proof of the first equality of condition (\ref{antid}) yields that for any 
$\theta\notin{\moic}:={\lcnsq}\cup{\qliani}$, where $\qliani$ is the null set outside of which the first equality of condition (\ref{antid}) holds true, the process 
$\Vp(\theta)$ 
is a 
martingale in $\mathcal{L}^1(Q_{\theta})$, since $\E_{Q_{\theta}}[S_t]=tg(\theta)\E_{Q_{\theta}}[X_1]=tg(\theta)\E_Q[X_1]<\infty$ for each $t\in\R_+$, where the 
second equality is an immediate consequence of the first equality of condition (\ref{antid}); hence $(i)$ follows for $\ell=1$. 
For $\ell=2$ applying \cite{sch}, Corollary 5.2.11, we get 
$S_t\in\mathcal{L}^2(Q_{\theta})$; 
hence $\Vp_t(\theta)\in\mathcal{L}^2(Q_{\theta})$ for each $t\in\R_+$. 

Ad $(ii)$: Let $Q\in\mathcal{M}^{*,\ell}_{\sagp,g}$. Applying Wald's identities (cf. e.g. \cite{sch}, Lemma 5.2.10) we get 
$S_t\in\mathcal{L}^{\ell}(Q)$; hence $\Vp_t(\vT)\in\mathcal{L}^{\ell}(Q)$ for any $t\in\R_+$.

Since $\agp$ has $Q_{\theta}$-stationary independent increments for any $\theta\notin{\moic}$, it follows by 
\cite{lm1}, Corollary 4.2 together with \cite{lm1err}
that it has $Q$-conditionally stationary independent increments, implying for all $s,t\in\R_+$ with $s\leq t$ that $S_t-S_s$ is $Q$-conditionally independent of 
$\haf_s^{\sagp}$; 
hence of ${\haf}_s$ (see \cite{lm1}, Lemma 4.7). The latter together with \cite{ct}, Section 7.3, Theorem 1 implies that condition $\int_A\E_Q[S_t-S_s\mid\haf_s]dQ=\int_A\E_Q[S_t-S_s\mid\vT]dQ$ holds true for all $s,t\in\R_+$ with $s\leq t$ and $A\in\haf_s$; hence inclusion $\mathcal{M}^{*,\ell}_{\sagp,g}\subseteq\M^{\ell}_{S,g}(\Vp(\vT))$ 
follows.

Conversely, if $Q\in\M^{\ell}_{S,g}(\Vp(\vT))$ 
then the process $\Vp(\vT)$ is a 
martingale in $\mathcal{L}^{\ell}(Q)$; hence $Q\in\mathcal{M}_{\sagp,g}^{*,\ell}$. 

Ad $(iii)$: 
Let $Q\in\mathcal{M}^{*,1}_{\sagp,g}$. 
If $\Vp$ 
is a 
martingale in $\mathcal{L}^1(Q)$, 
we then get that $\int_D(S_t-S_s)dQ=\int_D\E_Q[(S_t-S_s)]dQ$ for \mcol{all} $s,t\in\R_+$ with $s\leq t$ and for each $D\in\sigma(\vT)$. The latter together with the fact that by assumption $\agp$ is a $Q$-CMPP$(g(\vT),Q_{X_1})$ such that $g(\vT)$ and $X_1$ are $Q$-integrable, yields for each $t\in\R_+$ that
$$
\E_Q[S_t\mid\vT]={\E}_Q[S_t]\Longleftrightarrow 
tg(\vT)\E_Q[X_1]=t\E_Q[g(\vT)]\E_Q[X_1]\Longleftrightarrow g(\vT)=\E_Q[g(\vT)],
$$
where all equalities hold true $Q\upharpoonright\sigma(\vT)$-a.s.. Consequently, there exists a $\theta_0\in\vY$ such that 
$1=Q_{g(\vT)}(\{g(\theta_0)\})=Q\bigl(\vT^{-1}(g^{-1}(\{g(\theta_0)\})\bigr)
=Q_{\vT}(\{\theta:g(\theta)=g(\theta_0)\})$; hence such that $P_{\vT}(\{\theta:g(\theta)=g(\theta_0)\})=1$, since $Q_{\vT}\sim P_{\vT}$. 
Thus, $\agp$ is a $Q$-CPP$(g(\theta_0),Q_{X_1})$.
For $Q\in\mathcal{M}^{*,2}_{\sagp,g}$ if $\Vp$ is a 
martingale in $\mathcal{L}^2(Q)$ then $Q\in\mathcal{M}^{*,1}_{\sagp,g}$; 
hence $g(\vT)$ is degenerate at $g(\theta_0)$ for some $\theta_0>0$ as above.
The inverse implication is clear.\hfill$\Box$

\begin{prop}\label{cqr1}
Let $(\vO,\vS,P)$ be an arbitrary probability  
space and let 
$\sagp$ be an arbitrary aggregate process induced by a counting $\cnp$ and a size process 
$\csp$ such that 
$X_1^{\ell}\in\mathcal{L}^1(P)$ 
for $\ell\in\{1,2\}$. Consider the following assertions: 
\begin{enumerate}
\item 
The process 
$Y(\vT):=\{Y_t(\vT)\}_{t\in\R_+}:=\{\sagp_t-t\vT{\E_P[X_1\mid\vT]}\}_{t\in\R_+}$ 
is a {$P$-}martingale in $\mathcal{L}^{\ell}(P)$; 
\item
there exists a $P_{\vT}$-null set ${K}\in\mf{B}(\vY)$ such that for any $\theta\notin{K}$ the process 
$Y(\theta):=\{Y_t(\theta)\}_{t\in\R_+}:=\{\sagp_t-t\theta\E_{P_{\theta}}[X_1]\}_{t\in\R_+}$ is a 
{$P_{\theta}$-}martingale in $\mathcal{L}^{\ell}(P_{\theta})$. 
\end{enumerate}
Then $(i)\Longrightarrow(ii)$. If in addition $Y_t(\vT)^{\ell}\in\mathcal{L}^{1}(P)$ 
for any $t\in\R_+$ then assertions $(i)$ and $(ii)$ are equivalent.
\end{prop}

{\bf Proof.} 
Clearly, both processes in items $(i)$ and $(ii)$ are adapted to the filtration $\mathcal{F}$. Let us fix on arbitrary $\ell\in\{1,2\}$.
\smallskip

Ad $(i)\Longrightarrow(ii)$: 
We divide the proof in the following three steps.

{\bf (a)} 
There exists a $P_{\vT}$-null set ${{K}}^{\prime}\in\mf{B}(\vY)$ such that 
$Y_t(\theta)\in\mathcal{L}^{\ell}(P_{\theta})$ for every $\theta\notin{{K}}^{\prime}$ and $t\in\R_+$.

In fact, for given $t\in\R_+$ it follows by assumption $Y_t(\vT)^{\ell}\in\mathcal{L}^1(P)$ that 
$\sagp^{\ell}_t\in\mathcal{L}^1(P)$;
hence
$$
\infty>\int_{\vT^{-1}(B)}\sagp^{\ell}_tdP=\int_{\vT^{-1}(B)}\E_P[\sagp^{\ell}_t\mid\vT]dP=\int_B\E_{P_{\theta}}[\sagp^{\ell}_t]P_{\vT}(d\theta)
$$
for all $B\in\mf{B}(\vY)$, where the second equality follows by \cite{lm1}, Lemma 3.5, $(i)$. 
Consequently, there exists a 
{$P_{\vT}$}-null set ${K}_t\in\mf{B}(\vY)$ such that $\E_{P_{\theta}}[\sagp^{\ell}_t]<\infty$ for any $\theta\notin{K}_t$. 
Writing ${K}_1:=\bigcup_{t\in\mathbb{Q}_+}{K}_t$ we obtain $P_{\vT}({K_1})=0$ and 
$\E_{P_{\theta}}[\sagp^{\ell}_u]<\infty$ for any $\theta\notin{{K}}_1$ and $u\in\mathbb{Q}_+$. 

We then get for any given $t\in\R_+$ and $\theta\notin{K_1}$ that $\E_{P_{\theta}}[\sagp^{\ell}_t]<\infty$, since there exists a number 
$u\in\mathbb{Q}_+$ with $t\leq{u}$, implying by the monotonicity of the paths of $\sagp$ that $\sagp^{\ell}_t\leq\sagp^{\ell}_u$; 
hence $\E_{P_{\theta}}[\sagp^{\ell}_t]\leq\E_{P_{\theta}}[\sagp^{\ell}_u]<\infty$. 
  
But since by assumption $X_1^{\ell}\in\mathcal{L}^1(P)$, there exists a $P_{\vT}$-null set ${K}_2\in\mf{B}(\vY)$ such that 
$X_1^{\ell}\in\mathcal{L}^1(P_{\theta})$ for all $\theta\notin{K}_2$, writing ${K}^{\prime}:={K}_1\cup{K}_2$ we obtain that 
$Y_t(\theta)^{\ell}\in\mathcal{L}^1(P_{\theta})$ for any $\theta\notin{K}^{\prime}$.

{\bf (b)}
Let us fix arbitrary $u,t\in\R_+$ with $u\leq t$ and $A\in\mathcal {F}_u$. Define the real valued function $g_t$ on $\vO$ by means of $g_t(\omega):=f_t\circ(\mbox{id}_{\vO}\otimes\vT)(\omega)$ for each $\omega\in\vO$, where $f_t$ is the real valued function on $\vO\times\R$ given by $f_t(\omega,\theta):=\chi_A(\omega)Y_t(\theta)(\omega)$ for each $(\omega,\theta)\in\vO\times\R$. Also denote by $M$ the image measure of $P$ under $\mbox{id}_{\vO}\otimes\vT$.
Then the function $f_t^{\ell}$ is $M$-integrable, since 
$\int|f_t|dM=\int|g_t|dP=\int\chi_A|Y_t(\vT)|^{\ell}dP \leq \E_P[|Y_t{(\vT)}|^{\ell}]<\infty$,
where the first equality follows from \cite{lm1}, Proposition 3.8, $(ii)$.
\smallskip

{\bf (c)}
There exists a $P_{\vT}$-null set ${K}^{\prime \prime}\in\mf{B}(\vY)$ 
such that for any $\theta\notin{K}^{\prime \prime}$ the process 
{$Y(\theta)$} satisfies the martingale property.

In fact, let us fix on arbitrary $u,t\in\R_+$ with $u\leq{t}$, $A\in\mathcal{F}_u$ and $B\in\mf{B}(\vY)$. 
Writing $\wt{f}_t(\omega,\theta):=\chi_B(\theta)f_t(\omega,\theta)$ for each $(\omega,\theta)\in\vO\times\R$ and 
$\wt{g}_t:=\wt{f}_t\circ(\mbox{id}_{\vO}\otimes\vT)$ and applying (b) we get $\wt{f}_t\in{\mathcal L}^1(M)$. 
Therefore, we may apply \cite{lm1}, Proposition 3.8, $(ii)$, to infer

\begin{equation}
\int_B\int{f}_t^{\theta}dP_{\theta}P_{\vT}(d\theta)=\int_{\vT^{-1}(B)}g_tdP.
\label{x22a}	
\end{equation}

Applying again \cite{lm1}, Lemma 3.5, $(i)$ 
we obtain 
\begin{eqnarray*}
\lefteqn{{Y(\vT)}\quad\mbox{satisfies the martingale property}}\\
&\Longleftrightarrow&
\int_{A}Y_t(\vT)dP=\int_{A}Y_{u}(\vT)dP
\Longleftrightarrow
\int_{\vT^{-1}(B)}g_tdP=\int_{\vT^{-1}(B)}g_udP\\
&\stackrel{(\ref{x22a})}{\Longleftrightarrow}&
\int_B\int{f}_t^{\theta}
dP_{\theta}P_{\vT}(d\theta)
=\int_B\int{f}_u^{\theta}dP_{\theta}P_{\vT}(d\theta)\\
&\Longleftrightarrow&
\int_B\int_AY_t(\theta)dP_{\theta}{P_{\vT}}(d\theta)=\int_B\int_AY_u(\theta)dP_{\theta}{P_{\vT}}(d\theta)\adc{,}
\end{eqnarray*}
where all equivalences hold true for all $A\in\mathcal{F}_s$ 
and $B\in{\mf B}$, as both of them are arbitrary.
But the last equality is equivalent to the fact that  
for all $u,t\in\R_+$ with $u\leq t$ and for each $A\in{{\mathcal{F}}_u}$ such that
\begin{equation}\label{sm}
{\E_{P_{\theta}}[\chi_A(S_t-S_u)]=(t-u)\theta\E_{P_{\theta}}[X_1]P_{\theta}(A)}
\end{equation}
for any $\theta\notin{{K}}^{\prime \prime}_{A,u,t}$.

Denoting by $\mathcal{G}_u$ a countable algebra generating 
${{\mathcal{F}}_u}$, putting 
$$
K^{\prime \prime}:=\bigcup_{A\in\mathcal{G}_u}\bigcup_{s,t\in\mathbb{Q}_+,\;s\leq{t}}K^{\prime \prime}_{A,u,t}
$$ 
and applying a monotone class argument, we get a $P_{\vT}$-null set in $\mf{B}(\vY)$ such that for any $\theta\notin{{K}}^{\prime \prime}$ 
condition (\ref{sm}) is valid for any $A\in{\mathcal{F}}_u$ and for all $u,t\in\mathbb{Q}_+$ with $u\leq{t}$.

Put $\mathcal{H}_t:=\bigcap_{t<u\in\R_+}\mathcal{F}_u$ for any $t\in\R_+$. 
Then if we take $u,t\in\R_+$ with $u\leq{t}$, and if we write (\ref{sm}) for $u^{\prime}, t^{\prime}\in\mathbb{Q}$ with 
$u^{\prime}>u$ and $t^{\prime}>t$, and then let $u^{\prime}\downarrow{u}$ and $t^{\prime}\downarrow{t}$, we see that (\ref{sm}) holds true 
for all $A\in\bigcap_{{u^{\prime}>u}}\mathcal{F}_{u^{\prime}}=\mathcal{H}_u$, 
where the equality is a consequence of the following

\begin{it}  
Claim. For any $t\in\R_+$ we have $\bigcap_{u^{\prime}\in\Q,\;u^{\prime}\downarrow{t}}\mathcal{F}_{u^{\prime}}=\mathcal{H}_t$.

Proof.
\end{it}
For an arbitrary but fixed $t\in\R_+$ we have
$$
\mathcal{H}_t\subseteq\bigcap_{u^{\prime}\in\Q,\;u^{\prime}>t}{\mathcal{F}_{u^{\prime}}}
\subseteq\bigcap_{u^{\prime}\in\Q,\;u^{\prime}\downarrow{t}}\mathcal{F}_{u^{\prime}}.
$$ 

Conversely, for arbitrary $A\in\bigcap_{u^{\prime}\in\Q,u^{\prime}\downarrow{t}}\mathcal{F}_{u^{\prime}}$ we get $A\in\mathcal{F}_{u^{\prime}}$ for all $u^{\prime}\in\Q$ with $u^{\prime}\downarrow{t}$. Let $u$ be an arbitrary element of $R_+$ with $t<u$. Since the set $\Q$ is dense in $\R$, there exist $u^{\prime}\in\Q$ such that $t<u^{\prime}<u$ and $u^{\prime}\downarrow{t}$; hence $A\in\mathcal{F}_u$ because $A\in\mathcal{F}_{u^{\prime}}\sq\mathcal{F}_u$. As $u$ and $A$ are arbitrary we obtain $A\in\mathcal{H}_t$. \hfill$\Box$

Consequence, $Y(\theta)$ satisfies condition (\ref{sm}) for any $A\in\mathcal{H}_u$; hence for any $A\in\mathcal{F}_u$ since $\mathcal{F}_u\sq\mathcal{H}_u$, and any $\theta\notin{{K}}^{\prime \prime}$. Thus,assertion $(ii)$ follows.

Note that  $\mathcal{H}:=\{\mathcal{H}_t\}_{t\in\R_+}$ is a right continuous filtration such that $\mathcal{H}_{\infty}:=\sigma(\mathcal{H})=\mathcal{F}_{\infty}=\vS$.

Ad $(ii)\Longrightarrow (i)$: Since 
$Y_t(\vT)^{\ell}\in\mathcal{L}^1(P)$ 
we are left to show that $Y(\vT)$ satisfies the martingale property. But since $Y(\theta)$ satisfies condition (\ref{sm}) for all $u,t\in\R_+$ with $u\leq{t}$, all $\theta\notin{K}^{\prime \prime}$ and all $A\in\mathcal{F}_u$ we infer that
$$
\int_B\int_AY_t(\theta)dP_{\theta}P_{\vT}(d\theta)=\int_B\int_AY_u(\theta)dP_{\theta}P_{\vT}(d\theta)
$$
holds for all $B\in\mf{B}(\vY)$ true,  
equivalently (see the proof of step (c)), 
that $Y(\vT)$ fulfills the martingale property. Thus, $(i)$ follows.\hfill$\Box$

The process $Y(\theta)$ of 
the above proposition is known in ruin theory as {\em claim surplus process}.
Accordingly, the process $\Vp_{\T}(\vT)$ appearing in the next result can be regarded
as a {\em canonical} one, since it corresponds to the claim surplus process {$\Vp_{\T}(\theta)$}.

\begin{thm}\label{qm00}
Let $P\in\mathcal{M}^{*,2}_{\sagp}$. 
For each pair $(\beta,\xif)\in\mathcal{F}_P^2\times{\mathcal{R}^{*,2}_+(\vY)}$ there exist
\begin{enumerate}
\item 
a unique pair $(g,Q)\in\mf{M}_+(\vY)\times\mathcal{M}^{*,2}_{\sagp,g}$ satisfying conditions $(**)$ and $(M_{\xif})$ and such that the 
process $\Vp_{\T}(\vT):=\{\Vp_t(\vT)\}_{t\in\T}$ satisfies the property (NFLVR);
\item
an essentially unique r.c.p. $\{Q_{\theta}\}_{\theta\in\vY}$ of $Q$ over $Q_{\vT}$ consistent with $\vT$ satisfying for any 
$\theta\notin{\wt{V}_{*}}:={\qans}\cup{\moic}\in\mf{B}(\vY)$ conditions $Q_{\theta}\in\wt{\mathcal{M}}^2_{\sagp,g}$ 
and $(M_{\theta})$ and such that the process $\Vp_{\T}(\theta):=\{\Vp_t(\theta)\}_{\theta\in\T}$ satisfies condition (NFLVR). 
\end{enumerate}
\end{thm}

{\bf Proof.} Let be given a pair $(\beta,\xif)\in\mathcal{F}_P^2\times{\mathcal{R}^{*,2}_+(\vY)}$.

Ad $(i)$: By Corollary \ref{qm0}, $(ii)$, there exists a unique pair $(g,Q)\in\mf{M}_+(\vY)\times\mathcal{M}^{*,2}_{\sagp,g}$ satisfying conditions 
$(*)$ and $(M_{\xif})$. Applying Proposition \ref{cqr}, $(ii)$, we then infer that for any $T>0$ the process $\Vp_{\T}(\vT)$ is a 
$\haf_{\T}$-martingale in $\mathcal{L}^2(Q_T)$; hence it is a $\haf_{\T}$-semi-martingale in $\mathcal{L}^2(Q_T)$
(cf. e.g. \cite{ww}, Definition 7.1.1). 
The latter implies that $\Vp_{\T}(\vT)$ is a $\haf_{\T}$-semi-martingale in
$\mathcal{L}^2({P_T})$
since $P_T\sim Q_T$ (see \cite{ww}, Theorem 10.1.8).  
Thus, we may apply the Fundamental Theorem of Asset Pricing of Delbaen-Schachermayer (cf. \cite{ds}, Theorem 14.1.1), 
to deduce that the process $\Vp_{\T}(\vT)$ satisfies condition (NFLVR). 
\smallskip

Ad $(ii)$: By Corollary \ref{qm0}, $(iii)$, there exists an essentially unique r.c.p.
of $Q$ over $Q_{\vT}$ consistent with $\vT$ satisfying for any $\theta\notin{\qans}$ conditions $Q_{\theta}\in\wt{\mathcal{M}}^2_{\sagp,g}$ and $(M_{\theta})$. But by Proposition \ref{cqr}, $(i)$, there exists a $Q_{\vT}$-null set ${\moic}\in\mf{B}(\vY)$ such that for any $\theta\notin{\moic}$ the process for any $\theta\notin{\moic}\in\mf{B}(\vY)$ the process 
$\Vp_{\T}(\theta)$ is a 
$\haf_{\T}$-martingale in $\mathcal{L}^2\bigl((Q_{\theta})_T\bigr)$. 
Consequently, applying the same arguments as in $(i)$ we get that 
assertion $(ii)$ holds for any $\theta\notin{\wt{V}_{*}}$ true.\hfill$\Box$

\begin{rem}\label{44r}
\normalfont
It is common place that 
a stochastic integral with respect to a semi-martingale is defined 
under P.A. Meyer's usual conditions (cf. e.g. \cite{ww}, Definition 2.1.5). 
These conditions are not necessary, though, for the definition of such an integral 
(see \cite{ww}, pages 22, 23 and 150). 
Thus, the aforementioned integrals can be defined without them; hence 
the implication (EMM)$\Longrightarrow$(NFLVR) of Delbaen-Schachermayer's Fundamental Theorem of Asset Pricing 
(which was exploited for proving Theorem \ref{qm00}) remains valid even if the usual conditions are violated.
\end{rem}

\section{Applications to premium calculation principles}\label{newfapp} 

In this section, we briefly discuss potential 
implications of our main result and its consequences 
(see Theorem \ref{qQ} and Section \ref{cmppm}, respectively) in 
the computation of premium calculation principles. 


In classical Risk Theory a $P$-CPP $\sagp_{\T}:=\{S_t\}_{t\in[0,T]}$ is regarded as a proper model for describing 
the real or subjective behavior of the total risk of a period $\T:=[0,T]$, $T>0$, undertaken by an insurance company. According to the {\em financial pricing of insurance} (FPI for short) approach, introduced by Delbaen \& Haezendonck in \cite{dh}, at each time $t$ 
the premium $p_t$ assigned to the remaining risk $S_{T}-S_t$ of the period $(t,T]$ is given by means of 
$p_t=(T-t)p(Q)<\infty$ for each $t\in\T$, where $Q$ is a {\bf risk-neutral measure} (i.e. it is equivalent at each time moment $t$ with the 
original probability measure $P$ and such that the underlying price process becomes under $Q$ a martingale) such that
$p(Q)={\E}_Q[S_1]<\infty$, where $p(Q)$ denotes the monetary payout per time unit for holding a risk, that is the {\bf premium density}. 
Furthermore, $Q$ should also satisfy condition 
$p(P)<p(Q)<\infty$, equivalently 
\begin{equation}\label{ccc}
\E_P[S_t]<\E_Q[S_t]<\infty\quad\mbox{for all}\quad t>0.
\end{equation} 
For more details, see \cite{dh}, pages 269-271.

Under the FPI framework, 
we say that a {\bf premium calculation principle} (PCP for short) is a probability measure $Q$ on $(\vO,\vS)$ such that $Q\stackrel{pr}{\sim}P$, the family $\agp$ is a $Q$-CPP and $X_1\in\mathcal{L}^{1}(Q)$.

Note that in comparison with \cite{dh}, Definition 3.1, here the class of all PCPs is enlarged, since 
$\vS={\haf}_{\infty}\supseteq\haf_{\infty}^{\sagp}$. 
Theorem \ref{qQ} and its proof together with Remark \ref{qr} implies that 
for $P_{\vT}$-a.a. $\theta\in\vY$ the probability measures $Q_{\theta}$ are PCPs. 
Thus, every $Q\in{\mathcal{M}_{\sagp,g}^{1}}$ 
is a mixture of PCPs, while it can be easily seen that it is not a PCP itself. 

\smallskip

\begin{ex}\label{ex0} 
\normalfont
Let $P\in\mathcal{M}_S^{*,1}$. {\bf(a)} If $\beta(x,\theta)=cx-\ln\E_P[e^{cX}]$ for each $x,\theta\in\vY$ with $c>0$ 
then $\E_P[e^{cX_1}]<\infty$ and $\E_P[X_1e^{cX_1}]<\infty$; hence $\beta\in\mathcal{F}^1_P$ with $\alpha=0$ and so $g=id_{\vY}$.

Then for any $\xif\in\mathcal{R}_{+}^{*,1}(\vY)$ we get a
measure $Q\in\mathcal{M}^{*,1}_{\sagp}$ by Corollary \ref{qm0} such that the corresponding measure $Q_{\theta}$ is an Esscher PCP for all $\theta\notin\qans$, where $\qans$ is the $P_{\vT}$-null set of Corollary \ref{qm0}.
Since the function $\beta$ consists exclusively of the summand $\gamma$, the computation of $p(Q_{\theta})$, $\theta\notin\qans$,
will require nothing more in practice than repeating the computations made in Example 3.3 of \cite{dh}. 
Therefore, since $\E_Q[X_1]<\infty$ by $\beta\in\mathcal{F}_P^1$, it can be easily seen that 
for any $\theta\notin\qans$
\begin{equation}\label{ccc2}
p(P_{\theta})<p(Q_{\theta})<\infty
\end{equation}
if and only if 
$\E_P[X_1]\E_P[e^{cX_1}]<\E_P[X_1e^{cX_1}]$.
Since $\E_P[\vT]<\infty$ by $P\in\mathcal{M}_{\sagp}^{*,1}$, condition (\ref{ccc}) 
holds true if and only if 
$\E_P[\vT]\E_P[X_1]\E_P[e^{cX_1}]<\E_P[\vT\xif(\vT)]\E_P[X_1e^{cX_1}]$.\\
{\bf(b)} Likewise, if 
$\beta=c\in\R$ and $\gamma=0$ then $g=e^c id_{\vY}$ and for every $\xif\in\mathcal{R}_{+}^{*,1}(\vY)$ 
we obtain a measure $Q\in\mathcal{M}^{*,1}_{\sagp,g}$ by 
 Corollary \ref{qm0} such that the corresponding measure $Q_{\theta}$ is an Expected Value PCP for all $\theta\notin\qans$, while condition (\ref{ccc2}) holds if and only if $c>0$.
\end{ex}
In the above example we rediscovered both the Esscher's and the Expected Value PCPs, and 
extracted, for each case, not just a single PCP but a whole family $\{Q_{\theta}\}_{\theta\notin\qans}$ of such PCPs. 
Next we present some PCPs induced by pairs $(\beta,\xif)\in\mathcal{F}_P^2\times\mathcal{R}_{+}^{*,2}(\vY)$ such that the function 
$\beta$ consists of two non zero summands $\alpha$ and $\gamma$. 
To this aim recall that $\mathbf{Ga}(\mathrm{a,b})$, $\mathbf{Be}(\mathrm{a,b})$ and $\mathbf{U}(\mathrm{c,d})$ 
denote the law of Gamma, Beta and Uniform distribution with parameters $\mathrm{a,b}>0$ and $\mathrm{c,d}\in\R$ (cf. e.g. \cite{sch}). 
\begin{ex}\label{my1n} 
\normalfont
Let $P\in\mathcal{M}_S^{*,2}$. If $\beta(x,\theta)=\ln{\frac{x\theta}{\E_P[{X_1}]}}$ for each $x,\theta\in\vY$ then $\gamma(x)={-\ln\frac{\E_P[{X_1}]}{x}}$ and $\alpha(\theta)=\ln\theta$ for each $x,\theta\in\vY$; hence $\beta\in\mathcal{F}_P$ .
If $P_{X_1}=\mathbf{Exp}(1/5)$, $P_{\vT}=\mathbf{Ga}(2,2)$ and $\xif(\theta)=(27/8)\theta^2e^{-\theta}$ for $\theta\in\vY$ then 
$\E_P[X_1]=5$, $\E_P[X_1^2]=50$, $\E_P[X_1^3]=750$, and $\xif\in\mathcal{R}_+(\vY)$. 
Thus, by Theorem \ref{qQ} there exists
a pair $(g,Q)\in\mathfrak{M}_+(\vY)\times\mathcal{M}_{\sagp,g}$ such that 
$g(\theta)=\theta^2$ for all $\theta\in\vY$ and $Q_{\vT}=\mathbf{Ga}(3,4)$. 
Moreover, $\E_P[X_1^2e^{\gamma(X_1)}]=\E_P[X_1^3]/\E_P[X_1]=75$ and 
${\E_P[g^2(\vT)\xif(\vT)]}=\Gamma(8)/(3^5\cdot2)=840/81$; hence $(\beta,\xif)\in\mathcal{F}_P^2\times\mathcal{R}_{+}^{*,2}(\vY)$. 
It then follows by Corollary \ref{qm0} that the pair $(g,Q)\in\mathfrak{M}_+(\vY)\times\mathcal{M}^{*,2}_{\sagp}$ such that
for all $\theta\notin\qans$ the corresponding measure $Q_{\theta}$ is a PCP  
since 
$\E_{Q_{\theta}}[X_1]=\E_{P_{\theta}}[X_1e^{\gamma(X_1)}]={\E_{P_{\theta}}[{X_1^2}]/\E_{P_{\theta}}[{X_1}]}=
{\E_P[{X_1^2}]/\E_P[{X_1}]}=10$.
In addition, we have 
$\E_{Q_{\theta}}[N_1]=\theta^2=\theta\E_{P_{\theta}}[N_1]$ and so 
$p(Q_{\theta})=\E_{Q_{\theta}}[S_1]=\E_{Q_{\theta}}[N_1]\E_{Q_{\theta}}[X_1]=10\theta^2<\infty$. 
Since $p(P_{\theta})=\E_{P_{\theta}}[S_1]=\E_{P_{\theta}}[N_1]\E_{P_{\theta}}[X_1]=\theta\E_{P}[X_1]=5\theta<\infty$, 
it follows that condition (\ref{ccc2}) holds true if and only if $\theta>1/2$.
Furthermore, $\E_{Q}[N_1]={\E_Q[\vT^2]=\E_P[\vT^2\xif(\vT)]}=(81/24)\Gamma(5)=81$ and
$\E_{Q}[X_1]={\E_{P}[{X_1^2}]/\E_P[{X_1}]}=10$, implying 
$p(Q)=\E_Q[S_1]=\E_Q[N_1]\E_Q[X_1]=810>5=1\cdot5=\E_P[\vT]\E_P[X_1]=\E_P[N_1]\E_P[X_1]=p(P)$, 
that is condition (\ref{ccc}). 
Finally, by Theorem \ref{qm00} we deduce that the process $\{\Vp_t(\vT)\}_{t\in\T}$, where
$\Vp_t(\vT)=S_t-tg(\vT)\E_P[X_1e^{\gamma(X_1)}]=S_t-10t\vT^2$, satisfies condition (NFLVR). 
\end{ex} 

\begin{ex}\label{my2n} 
\normalfont
Let $P\in\mathcal{M}_S^{*,2}$. If $\beta(x,\theta)=cx+\ln[(c+\theta)\E_{P_{\theta}}[e^{-\theta X_1}]/(c+1)^2]$ for each $x,\theta,c\in\vY$ such that $c+3>(c+2)^2\ln[(c+2)/(c+1)]$, then $\gamma(x)=cx-2\ln(c+1)$ and $\alpha(\theta)=\ln(c+\theta)+\ln\E_{P_{\theta}}[e^{-\theta X_1}]$  for each $x,\theta,c$ as above.
If $P_{X_1}=\mathbf{Ga}(c+1,2)$, $P_{\vT}=\mathbf{Be}(2,1)$ and $\xif(\theta)=1/2\theta$ for each $\theta\in(0,1)$, then 
$\E_P[X_1^2]=6(c+1)^{-2}$, $\E_P[e^{cX_1}]=(c+1)^2$, $\E_P[e^{\gamma(X_1)}]=(c+1)^{-2}\E_P[e^{cX_1}]=1=\E_P[\xif(\vT)]$; 
hence $(\beta,\xif)\in\mathcal{F}_P\times\mathcal{R}_+(\vY)$. 
Thus, by Theorem \ref{qQ} there exists a pair $(g,Q)\in\mathfrak{M}_+(\vY)\times\mathcal{M}_{\sagp,g}$ such that 
$g(\theta)=(c+\theta)\theta\E_{P_{\theta}}[e^{-\theta X_1}]=(c+1)^2(c+1+\theta)^{-2}(c+\theta)\theta$ for all $\theta\in\vY$ 
and $Q_{\vT}=\mathbf{U}(0,1)$. Moreover, 
$\E_P[X_1^2e^{\gamma(X_1)}]=(c+1)^{-2}\E_P[X_1^2e^{cX_1}]=6$ and 
$$\E_P[g^2(\vT)\xif(\vT)]=
\frac{1}{2}\int_{0}^{1}\theta\frac{(c+\theta)^2(c+1)^4}{(c+\theta+1)^4}P_{\vT}(d\theta)
<\frac{(c+1)^4}{2},
$$
implying $(\beta,\xif)\in\mathcal{F}_P^2\times\mathcal{R}_{+}^{*,2}(\vY)$.
We then get by Corollary \ref{qm0} that the pair $(g,Q)\in\mathfrak{M}_+(\vY)\times\mathcal{M}^{*,2}_{\sagp}$ such that
for all $\theta\notin\qans$ the corresponding measure $Q_{\theta}$ is a PCP  
since 
$\E_{Q_{\theta}}[X_1]=\E_{P_{\theta}}[X_1e^{\gamma(X_1)}]=(c+1)^{-2}\E_P[X_1e^{cX_1}]=2$. In addition, 
$\E_{Q_{\theta}}[N_1]=g(\theta)$
$\E_{Q_{\theta}}[N_1]=g(\theta)$
 and so $p(Q_{\theta})=2(c+\theta)\theta[(c+1)^2/(c+1+\theta)^2]<\infty$. 
Since $p(P_{\theta})=\theta\E_{P}[X_1]=2\theta/(c+1)<\infty$, 
it follows that condition (\ref{ccc2}) holds true if and only if 
$$
\frac{2\theta}{c+1}<\frac{2(c+\theta)(c+1)^2\theta}{(c+\theta+1)^2}
\Longleftrightarrow\theta^2-(c+1)(c^2+2c-1)\theta+(c+1)^2(1-c-c^2)<0.
$$
Furthermore, $\E_{Q}[X_1]=\E_{Q_{\theta}}[X_1]=2$ for any $\theta\notin\qans$ and 
$\E_Q[N_1]=\E_Q[g(\vT)]=\int g(\theta)Q_{\vT}(d\theta)=(c+1)^2J(c)$, 
where 
$$
J(c):=\int_{0}^{1}\theta[(c+\theta)/(c+1+\theta)^2]d\theta=\frac{c+3}{c+2}+(c+2)\ln\frac{c+1}{c+2}\in(0,\infty);
$$
hence $p(Q)=2(c+1)^2J(c)<\infty$. 
Since it follows by our assumptions that $\E_P[\vT]=2/3$ and $\E_P[X_1]=2(c+1)^{-1}$ we get
$p(P)=(4/3)(c+1)^{-1}<\infty$. 
Thus condition (\ref{ccc}) holds true if and only if $(4/3)(c+1)^{-1}<2(c+1)^2J(c)$ or $J(c)<(2/3)(c+1)^{-3}$.

Finally, by Theorem \ref{qm00} we deduce that the process $\{\Vp_t(\vT)\}_{t\in\T}$, where
$\Vp_t(\vT)
=S_t-2(c+\vT)[(c+1)^2t\vT/(c+1+\vT)^2]$,
satisfies condition (NFLVR). 
\end{ex}
Note that in Example \ref{my1n}, the change of measure from $P_{\theta}$ to $Q_{\theta}$ as well as from $P$ to $Q$ has 
a multiplicative effect both to the expected number of claims per time unit (namely $\E_{P_{\theta}}[N_1]$, $\E_P[N_1]$) 
and the expected claim size. The latter effect is 
confirmed in Example \ref{my2n} for the claim size but a similar comment cannot be made for the expected number of claims 
as, in that case, the change of measure effect depends on the parameter $c$ of the claim size distribution 
$P_{X_1}$.

{\small
\begin{tabular}{ll}
{\tt D.P. Lyberopoulos {\sf and} N.D. Macheras}&$\;$\\
{\sf Department of Statistics and Insurance Science}&$\;$\\
{\sf University of Piraeus, 80 Karaoli and Dimitriou street}&$\;$\\
{\sf 185 34 Piraeus, Greece}&$\;$\\
{\sf E-mail:} {\tt dilyber@webmail.unipi.gr}$\;$ {\sf and}$\;$ {\tt macheras@unipi.gr}
\end{tabular}
}

\bigskip

\appendix
\section*{{Appendix: List of Null Sets}}

Next we present, for completeness and convenience, a list of the null sets met in the main body of the paper and used in the proofs of our results. This list is divided into three different tables and each table consists of three columns, indicating 
{\em the null set we refer to, the position in the text where the null set is met for the first time} and {\em the null set's property}. 
Unless it is stated otherwise, the property we refer to at each line is valid under the disintegrating measure $P_{\theta}$ and for any $\theta$ outside the corresponding null set.
Note also that there are null sets appearing only once in the text for technical reasons; these sets are given in brackets where thought necessary. 

\begin{center}
{\bf\large Section \ref{lm33}}
\end{center}
\begin{tabular}[t]{lcl}
${\lpns}$&\ref{11}, $(i)$& $\cnp$ independent of $\csp$\\
${\ins}$&{proof of \ref{11}, $(ii)$}&independent $\csp$\\
${\idnso}$&{proof of \ref{11}, $(ii)$}&identically distributed $\csp$\\
${\ldns}$&proof of \ref{11}, $(ii)$&i.i.d. $\csp$ ($\ldns:={\ins}\cup{\idnso}$)\\
${\lrns}$&proof of \ref{11}, $(iii)$&$\rp$ \rpw ($\lrns=\lpns\cup\ldns$)\\ 
${\lcns}$&proof of \ref{12}&$\agp$ CPP($\theta,(P_{\theta})_{X_1}$) ($\lcns=\lrns\cup\poins$)
\end{tabular}
 
\begin{center}
{\bf\large Section \ref{lm34}} 
\end{center}
\begin{tabular}[t]{lcl}
${\lcnsq}$&proof of \ref{tsp}, step (a)&$\agp$ $Q_{\theta}$-CPP$(g(\theta),(Q_{\theta})_{X_1})$\\
${\ccns}$&proof of \ref{tsp}, step (a)&${\ccns}:={\lcns}\cup{\lcnsq}$\\
${\ldnsq}$&proof of \ref{tsp}, step (b)&$Q_{\theta}$-i.i.d. $\csp$\\
${\qliani}$&proof of \ref{tsp}, step (b)&first equality of (\ref{antid})\\
$\;$&$\;$&{[}\textsl{${\machns}$ are also involved in the proof}{]}\\ 
${\gplians}$&proof of \ref{tsp}, step (b)&last equality of (\ref{antid})\\
$\;$&$\;$&{[}\textsl{${\umachns}$ are also involved in the proof}{]}
\end{tabular}

\begin{tabular}[t]{lcl}
${\lpnsq}$&proof of \ref{tsp}, step (b)& {condition (\ref{antid}) (${\caggns}:={\ldnsq}\cup{\qliani}\cup{\gplians}$)}\\
$\wt{G}_d$&proof of \ref{tsp}, step (d) & $\cnp$ $Q_{\theta}$-independent of $\csp$\\
\textsl{${\pliani}$}&proof of \ref{tsp}, step (d)&$\pliani=\ccns\cup\lpnsq\cup\wt{G}_d$\\
${\qcns}$&proof of \ref{tsp}, step (e)&$Q_{\theta}(\vT^{-1}(B))=\chi_{B}(\theta)$ 
{for each $B\in\mf{B}(\vY)$}\\
$\;$&$\;$&{(}$\{Q_{\theta}\}_{\theta{\in\vY}}$ consistent with $\vT${)}\\
$\pcns$&proof of \ref{tsp}, step (e)&$P_{\theta}(\vT^{-1}(B))=\chi_{B}(\theta)$ 
{for each $B\in\mf{B}(\vY)$}\\
$\;$&$\;$&{(}$\{P_{\theta}\}_{\theta{\in\vY}}$ consistent with $\vT${)}\\
${\qans}$&proof of \ref{tsp}, step (e)&${\qans}:={\pliani}\cup{\pcns}\cup{\qcns}$\\
$\wt{H}$&proof of \ref{46}&$\agp$ $Q_{\theta}$-CPP($g(\theta),(Q_{\theta})_{X_1}$), $P_{\theta}$-CPP{($\theta,(P_{\theta})_{X_1}$)}\\
$\wh{H}$&proof of \ref{46}&$P_{\theta}\neq Q_{\theta}$
\end{tabular}

\begin{center}
{\bf\large Sections \ref{lm35} and \ref{cmppm}}
\end{center}
\begin{tabular}[t]{lcl}
${\qans}$&proof of \ref{qQ}, $(ii)$& condition $(M_{\theta})$\\
$\wh{L}$&\ref{qr}, (d)&$\agp$ is a $Q_{\theta}$-CPP$(g(\theta),(Q_{\theta})_{X_1})$\\
$\moic$&proof of \ref{cqr}, $(i)$&{$\E_{Q_{\theta}}[X_1]=\E_Q[X_1]$ (${\moic}:={\lcnsq}\cup{\qliani}$)}\\
{$K^{\prime}$}&{proof of \ref{cqr1}}&{$\mathcal{L}^{1}(P_{\theta})$-integrability of $Y(\theta)$}\\
$\;$&$\;$&{[\textsl{$K_t$, $K_{A,u,t}$ are also involved in the proof}]}\\
{$K^{\prime\prime}$}&{proof of \ref{cqr1}}&{$(P_{\theta},\mathcal{H})$-martingale property of $Y(\theta)$}\\
{$K$}&{proof of \ref{cqr1}}&{$P_{\theta}$-martingale property of $Y(\theta)$ ($K:=K^{\prime}\cup K^{\prime\prime}$)}\\
{$\wt{V}_{*}$}&{proof of \ref{qm00}}&{$Q_{\theta}\in\wt{\mathcal{M}}^2_{\sagp,g}$, 
$(M_{\theta})$ and (NFLVR) for $\Vp_{\T}(\theta)$ ($\wt{V}_{*}:={\qans}\cup{\moic}$)}
\end{tabular}


{\small
\begin{thebibliography}{99}
\bibitem{ba} 
\auth{Bauer, H.}: 
\artb{Probability Theory}, Walter de Gruyter, Berlin-- New York (1996). 


\bibitem{bdw} 
\auth{Boogaert, P.} and \auth{de Waegenaere, A.} (1990):
\artp{Simulation of ruin probabilities.}
\artj{Insurance Math. Econom.} {\bf 9}, 95-99.

\bibitem{ct} 
\auth{Chow, Y.S.} and \auth{Teicher, H.}:
\artb{Probability Theory}. Second Edition, Springer - Verlag, New York (1988).  

\bibitem{dh} 
\auth{Delbaen, F.} and \auth{Haezendonck, J.} (1989): 
\artp{A martingale approach to premium calculation principles in an arbitrage free market}. 
\artj{Insurance Math. Econom.} {\bf 8}, 269-277.

\bibitem{ds}  
\auth{Delbaen, F.} and \auth{Schachermayer, W.}: 
\artb{The Mathematics of Arbitrage}. Springer-Verlag, Berlin-Heidelberg (2006).

 

\bibitem{em00}  
\auth{Embrechts, P.} (2000): 
\artp{Actuarial versus financial pricing of insurance}. 
\artj{J Risk Finance} \artn{1}, 17--26.

\bibitem{em97}  
\auth{Embrechts, P.} and \auth{Meister, S.} (1997): 
\artp{Pricing Insurance Derivatives, the Case of CAT-Futures}. 
Proceedings of the 1995 Bowles Symposium on Securitization of Risk, 
George State University Atlanta, Society of Actuaries, Monograph M-FI97-1, 15--26.

\bibitem{fa} 
\auth{Faden, A.M.} (1985):
\artp{The existence of regular conditional probabilities: Necessary and sufficient conditions.} 
\artj{Ann. Probab.} \artn{13} (No. 1), 288--298. 


\bibitem{fr4} 
\auth{Fremlin, D.H.:} 
\artb{Measure Theory, Vol. 4, Topological Measure Spaces}.
Torres Fremlin, Colchester (2003).


\bibitem{ho} 
\auth{Holtan, J.} (2007):
\artp{Pragmatic insurance option pricing.} 
\artj{Scand. Actuar. J.} \artn{2007} (1), 53--70. 



\bibitem{ks} 
\auth{Karatzas, I.} and \auth{Shreve, S.E.}: 
\artb{Brownian Motion and Stochastic Calculus}. Springer-Verlag, New York (1988). 

\bibitem{lm1} 
\auth{Lyberopoulos, D.P.} and \auth{Macheras, N.D.} (2012):
\artp{Some characterizations of mixed Poisson processes}. 
\artj{Sankhy\=a} \artn{74-A}, Part 1, 57--79.

\bibitem{lm5jmaslo} 
\auth{Lyberopoulos, D.P.} and \auth{Macheras, N.D.} (2013): 
\artp{A construction of mixed Poisson processes via disintegrations}. 
\artj{Math. Slovaca} \artn{63}, No. 1, 167--182.

\bibitem{lm1err} \auth{Lyberopoulos, D.P.} and \auth{Macheras, N.D.} (2014):
Erratum-Some characterizations of mixed Poisson processes. 
\artj{Sankhy\=a} \artn{76-A}, Part 1, 177.



\bibitem{mei}  
\auth{Meister, S.}: 
Contributions to the mathematics of catastrophe insurance futures.
Diplomarbeit, ETH-Z{\"u}rich (1995).

\bibitem{mol}  
\auth{M{\o}ller, T.} (2004): 
\artp{Stochastic orders in dynamic reinsurance markets}.
\artj{Finance Stochast.} \artn{8}, 479-499.

\bibitem{pa} 
\auth{Pachl, J.K.} (1978): 
\artp{Disintegration and compact measures}. 
\artj{Math. Scand.} \artn{43}, 157-168.



\bibitem{sch} \auth{Schmidt, K.D.}: 
\artb{Lectures on Risk Theory}. B.G. Teubner, Stuttgart (1996).


\bibitem{mms3} \auth{Strauss, W., Macheras, N.D.} and \auth{Musial, K.} (2004):
\artp{Splitting of liftings in products of probability spaces.} 
\artj{Ann. Probab.} \artn{32} (No. 3B), 2389--2408.


\bibitem{ww} \auth{von Weizs\"{a}cker} and \auth{H., Winkler, G.}: 
\artb{Stochastic Integrals. An Introduction}. Springer Fachmedien Wiesbaden (1990).

\bibitem{yual} 
\auth{Yu, S., Unger, A.J.A., Parker, B.} and \auth{Kim, T.} (2012):
\artp{Allocating risk capital for a brownfields redevelopment project under
hydrogeological and financial uncertainty.} 
\artj{J. Environ. Manag.} \artn{100}, 96--108. 


\end{thebibliography}}
\end{document}